\documentclass[11pt]{article}
\usepackage[utf8]{inputenc}

\usepackage{amsmath,amssymb}
\usepackage{amsthm}
\usepackage{enumitem}
\usepackage{appendix}
\usepackage{theoremref}
\usepackage{caption}
\usepackage{soul}

\usepackage[usenames]{color}

\usepackage[margin=1in]{geometry}
\usepackage{amsmath}
\usepackage{mathrsfs}
\usepackage{amsfonts}
\usepackage{cite}
\usepackage[mathscr]{euscript}
\usepackage{amssymb}
\usepackage{appendix}
\usepackage{tikz-cd}
\usepackage{empheq}
\usepackage{mathtools}
\usepackage{theoremref}
\usepackage[colorlinks=true, linkcolor=black, urlcolor=black, citecolor=black]{hyperref}
\usepackage{comment}
\usepackage{dsfont}
\usepackage{multicol}
\usepackage{graphicx}
\usepackage{float}
\usepackage{enumitem}
\usepackage{subcaption}
\usepackage{array}
\usepackage{mhsetup}
\usepackage[misc]{ifsym}
\usepackage{booktabs}

\usepackage{datetime}

\usepackage[in]{fullpage}
\usepackage{algorithm}
\usepackage{algorithmic}
\usepackage{authblk}

\graphicspath{{./ESFigures/}}

\newtheorem{thm}{Theorem}[section]

\newtheorem{defn}[thm]{Definition}

\newcommand{\R}{\mathbb{R}}

\newcommand{\N}{\mathbb{N}}
\newcommand{\Z}{\mathbb{Z}}

\newcommand{\z}{\mathbf{z}}  
\newcommand{\m}{\mathbf{m}}
\newcommand{\bc}{\mathbf{c}}
\newcommand{\bx}{{\mathbf x}}
\newcommand{\by}{{\mathbf y}}
\newcommand{\bu}{{\mathbf u}}
\newcommand{\bv}{\mathbf{v}}
\newcommand{\bl}{{\mathbf l}}
\newcommand{\bk}{\mathbf{k}}
\newcommand{\w}{\mathbf{w}}
\newcommand{\om}{\overline{m}}          

\newcommand{\obm}{\mathbf{\overline{m}}}
\newcommand{\obx}{\overline{\mathbf{x}}}
\newcommand{\oby}{\overline{\mathbf{y}}}
\newcommand{\obz}{\overline{\mathbf{z}}}

\newcommand{\rd}{\mathrm{d}}

\newcommand{\pdone}[2]{\dfrac{\partial {#1}}{\partial {#2}}}
\newcommand{\pdonet}[2]{\frac{\partial {#1}}{\partial {#2}}}

\newcommand{\per}{{\mathrm{per}}}
\newcommand{\omp}{{\omega'}}

\newcommand{\mres}{%
	\,\raisebox{-.127ex}{\reflectbox{\rotatebox[origin=br]{-90}{$\lnot$}}}\,%
}

\newcommand{\revA}[1]{\textcolor{black}{#1}}
\newcommand{\revB}[1]{\textcolor{black}{#1}}
\newcommand{\revC}[1]{\textcolor{black}{#1}}

\title{A new implementation of the geometric method\\ for solving the Eady slice equations}
\author[1]{C.~P.~Egan\footnote{\Letter\ Corresponding author. E-mail: cpe4@hw.ac.uk}}
\author[1]{D.~P.~Bourne}
\author[2]{C.~J.~Cotter}
\author[3]{M.~J.~P.~Cullen}
\author[1]{\\B.~Pelloni}
\author[4]{S.~M.~Roper}
\author[5]{M.~Wilkinson}
\affil[1]{Maxwell Institute for Mathematical Sciences and Department of Mathematics, Heriot-Watt University, Edinburgh, UK.}
\affil[2]{Department of Mathematics, Imperial College London, London, UK.}
\affil[3]{Met Office, Exeter, UK (retired)}
\affil[4]{School of Mathematics and Statistics,
University of Glasgow, Glasgow, UK.}
\affil[5]{Department of Mathematics, Nottingham Trent University, Nottingham, UK.}

\date{}

\begin{document}
\maketitle

\begin{abstract}
We present a new implementation of the geometric method of Cullen \& Purser (1984) for solving the semi-geostrophic Eady slice equations, which model large scale atmospheric flows and frontogenesis. The geometric method is a Lagrangian discretisation, where the PDE is approximated by a particle system. An important property of the discretisation is that it is energy conserving. We restate the geometric method in the language of semi-discrete optimal transport theory and exploit this to develop a fast implementation that combines the latest results from numerical optimal transport theory with a novel adaptive time-stepping scheme. Our results enable a controlled comparison between the Eady-Boussinesq vertical slice equations and their semi-geostrophic approximation. We provide further evidence that weak solutions of the Eady-Boussinesq vertical slice equations converge to weak solutions of the semi-geostrophic Eady slice equations as the Rossby number tends to zero.
\end{abstract}

\section{Introduction}

In this paper, we study a meshfree method for solving the \emph{semi-geostrophic Eady slice equations}, which are stated in equations \eqref{eqn:SGES1}, \eqref{eqn:vThRel0} in Eulerian coordinates, and in equation \eqref{eqn:SGESLagCts} in Lagrangian coordinates. The Eady slice model is a simplified model of large-scale (small Rossby number) atmospheric flow, and is capable of predicting the formation of atmospheric fronts \cite{cullen2021mathematics}. We develop an efficient implementation of the  \emph{geometric method} \cite{cullen1984extended} for solving this PDE using an adaptive time-stepping algorithm (Algorithm \ref{alg:adaptiveTimeStep}). As an application, we use it to test the validity of the semi-geostrophic approximation of the Euler-Boussinesq equations.

\paragraph{Background and motivation.}
It is important that numerical schemes used in atmospheric models adhere to the energy dynamics of the physical system that they represent. Otherwise, the variability of that system will not be accurately captured in the corresponding numerical solutions.
In Section \ref{sect:results}, we continue
the work of \cite{visram2014framework} and \cite{yamazaki2017vertical}, in which the idealised baroclinic lifecycle first introduced by Eady in \cite{eady1949long} is used as a test problem for numerical schemes for atmospheric models.
This problem has two particular features that make it a challenging and informative test. The first is that solutions form discontinuities, which are an idealised representation of atmospheric fronts. The second is that, in the absence of forcing, the total energy in the physical system is conserved in the large-scale limit.
An overview of existing numerical schemes for solving the Eady slice problem is given in \cite{visram2014framework}. The numerical scheme under consideration in this work is Cullen and Purser's \emph{geometric method} \cite{cullen1984extended} for solving the semi-geostrophic approximation of the Eady-Boussinesq vertical slice equations.

The geometric method was first introduced in \cite{cullen1984extended}.
It is a Lagrangian, meshfree method (particle method).
It is derived by first writing the semi-geostrophic Eady slice equations \eqref{eqn:SGES1},  \eqref{eqn:vThRel0} in Lagrangian coordinates (equation \eqref{eqn:SGESLagCts}), then by approximating
 the modified geopotential $P$  by a piecewise affine function (see equation \eqref{eqn:pwConstNablaP}). This reduces the PDE \eqref{eqn:SGESLagCts} to the system of ODEs given in equation \eqref{eqn:ODE}. The right-hand side of the ODE is defined in terms of a \emph{semi-discrete optimal transport problem}; to evaluate, it we must solve an optimal partitioning problem (Definition \ref{def:SDOT}) numerically or, equivalently, find the maximum of a concave function ($\mathcal{K}$, defined in equation \eqref{eqn:kappa}). We solve the ODE using a novel time-stepping algorithm (Algorithm \ref{alg:adaptiveTimeStep}), which is designed to 
accelerate the optimal transport solver.

{One of the main differences between our implementation of the geometric method and the original implementation in \cite{cullen1984extended} is that we restate it in the language of \emph{semi-discrete optimal transport theory} \cite[Chapter 4]{merigot2021optimal}, which allows us to exploit recent results from this field, such as the damped Newton method of Kitagawa, M\'erigot \& Thibert (2019) \cite{kitagawa2016convergence}. In our companion paper \cite{bourne2022semi}},
the pioneering design of the geometric method {is} put on a rigorous footing 
{for the closely-related three-dimensional semi-geostrophic equations}.
The {original} implementation of the geometric method in \cite{cullen1984extended} is now regarded as the first numerical method for solving the semi-discrete optimal transport problem \cite[Section 4.1]{merigot2021optimal}. The geometric method, however, was developed before the term \emph{semi-discrete optimal transport} was even coined. This paper sees a reversal of the knowledge exchange: we use the very latest results from semi-discrete optimal transport theory to improve the implementation of the geometric method. Our implementation builds on that of \cite{cullen1984extended} using 38 years of developments in computational geometry and optimal transport theory, including the fast algorithm and convergence guarantee of \cite{kitagawa2016convergence}. The robustness and conservation properties of the geometric method are optimised by our specific implementation choices.
Semi-discrete optimal transport has also been used to simulate the incompressible Euler equations \cite{gallouet2018lagrangian,levy2018notions,MerigotMirebeau2016}, barotropic fluids and porous media flow \cite{GallouetMerigotNatale2021}, and in astrophysical fluid dynamics \cite{LevyMohayaeevonHausegger2021}.

An advantage of the geometric method over finite element methods, such as those used in \cite{visram2014framework} and \cite{yamazaki2017vertical}, is that it is \emph{structure preserving}: 
{solutions of the discretised equations \eqref{eqn:ODE} conserve the total energy \eqref{eqn:totalGeoEnergy}, and give rise to mass-preserving flows in the fluid domain.}
Correspondingly, the numerical solutions that we obtain accurately conserve the total energy, even when frontal discontinuities form.
Similarly, numerical solutions obtained in \cite{culrou} by using the geometric method, and discussed further in \cite{cullen2007modelling}, exhibit non-dissipative lifecycles and show very little sensitivity to the discretisation, indicating high predictability. In contrast, numerical solutions of the Eulerian Eady slice equations obtained in \cite{visram2014framework, yamazaki2017vertical} by using Eulerian numerical methods are invariably dissipative. The diagnostics extracted in \cite[Section 3.3]{yamazaki2017vertical} show that they systematically violate energy conservation properties after front formation, resulting in a significant loss of total energy over multiple lifecycles.

The main application of our new algorithm is to obtain
numerical evidence that the semi-geostrophic Eady slice equations \eqref{eqn:SGES1}, \eqref{eqn:vThRel0}, in which geostrophic balance of the out-of-slice wind and hydrostatic balance are enforced, are the correct limit of the Eady-Boussinesq vertical slice equations \eqref{eqn:EadySlice1} {in the large-scale limit (i.e. as the Rossby number $\mathrm{Ro}$ tends to $0$)}. In the classical setting, solutions of the Eady-Boussinesq vertical slice equations with periodic boundary conditions converge {strongly in $L_2$} to solutions of the corresponding semi-geostrophic equations as $\mathrm{Ro}\to 0$ \cite[Theorem 4.5]{brenier13}. Correspondingly, the numerical solutions obtained in \cite{visram2014framework} satisfy geostrophic and hydrostatic balance, up to an error of order $\mathcal{O}(\mathrm{Ro}^2)$, while the solutions are smooth, and order $\mathcal{O}(\mathrm{Ro})$ after front formation. 
On the other hand, in both \cite{visram2014framework} and \cite{yamazaki2017vertical} it was remarked that there is a significant difference between the Eulerian solutions presented in those papers and the semi-geostrophic solution presented in \cite{culrou}, which cannot be accounted for by the dissipative nature of the former. Subsequent investigations revealed that the problems studied in the two papers are not identical. Our implementation of the geometric method enables a controlled comparison between the two approaches. 
Using the physical parameters from \cite{visram2014framework}, we obtain numerical solutions of the semi-geostrophic Eady slice equations
\eqref{eqn:SGES1}, \eqref{eqn:vThRel0}
and confirm that the differences between these and the Eulerian solutions obtained in \cite{visram2014framework} are consistent with the loss of Lagrangian conservation in the Eulerian solutions. This provides numerical evidence that the semi-geostrophic approximation is the correct small Rossby number approximation of the 
Eady-Boussinesq vertical
slice equations, and suggests that a result similar to \cite[Theorem 4.5]{brenier13} may also hold for weak solutions which, unlike classical solutions, can possess non-dissipative singularities that represent the evolution of weather fronts.

\paragraph{Outline of the paper.}
{In Section \ref{sec:Model}
we introduce the semi-geostrophic Eady slice equations 
in Eulerian coordinates (Section \ref{sec:Eulerian}) and Lagrangian coordinates (Section \ref{sec:Lagrangian}), as well as an important steady shear flow. In Section \ref{sect:geomMethod} we 
discretise the Eady slice equations using 
the geometric method. 
This leads us to the system of ODEs \eqref{eqn:ODE}. 
We give an algorithm for solving these ODEs in Section \ref{sect:implementation}, which involves discretising the initial data (Section \ref{subsec:quantizeIC}), evaluating the right-hand side of the ODE by solving a semi-discrete optimal transport problem (Section \ref{subsec:SolveOT}), and  using an adaptive time-stepping scheme (Section \ref{subsec:adaptive}). Section \ref{sect:results} includes some numerical experiments, where we study the stability of a steady shear flow. We illustrate front formation in Section \ref{Subsec:UNM} and study the validity of the semi-geostrophic approximation in Section \ref{subsec:ComparisonLiterature}. This numerical study is complemented by an analytical linear instability analysis in Appendix \ref{Sec:LinearInstabilityAnalysis}.
}

\section{Governing equations}
\label{sec:Model}

We study the semi-geostrophic approximation of the Eady-Boussinesq vertical slice model in an infinite channel of height $H$. We will consider solutions that are $2L$-periodic in the first coordinate direction. 
 We start by summarising the derivation of this model from the Euler-Boussinesq equations and we then draw comparisons with the Eady-Boussinesq vertical slice model considered in \cite{visram2014framework}.

\subsection{Eady slice model} 
\label{sec:Eulerian}
Our starting point is the Euler-Boussinesq equations for an incompressible fluid:
\begin{align}
\label{eqn:EulBous}
\pdone{\bu}{t}+\bu \cdot \nabla\bu + f \hat{\z} \times \bu &= -\nabla \Phi + \frac{g}{\theta_0}\Theta\hat{\z},
\\
\pdone{\Theta}{t}+\bu\cdot \nabla \Theta &= 0,
\\ 
\label{eqn:EulBous:3}
\nabla\cdot \bu &= 0.
\end{align}
{When applied to the atmosphere,} $\bu=(u,v,w)$ is the Eulerian fluid velocity, $\Theta$ is the potential temperature, $\Phi$ is the {geopotential}, $f$ is the Coriolis parameter, $g$ is the acceleration due to gravity, $\theta_0$ is a constant reference potential temperature, $\nabla=(\partial_x,\partial_y,\partial_z)$ is the gradient operator, and $\hat{\z}=(0,0,1)$. The three coordinates of a position vector $\bx=(x,y,z)\in\R^3$ represent longitude, latitude and altitude, respectively.

We decompose the potential temperature and 
{geopotential} as
\begin{align}
    \label{eqn:thDe}
	\Theta(\bx,t) & = \theta_0 + \bar{\theta}(y) + \theta(x,z,t),
\\
\label{eqn:phiDe}
	\Phi(\bx,t) & = \phi_0(z) + \bar{\phi}(y,z) + \phi(x,z,t), 
\end{align}
where the reference {geopotential} $\phi_0$ is given by 
\begin{align*}
	\phi_0(z)=gz,
\end{align*}
and the background potential temperature $\bar{\theta}$ and background {geopotential} $\bar{\phi}$ are given by
\begin{align*}
	\bar{\theta}(y)= s y,\qquad \bar{\phi}(y,z)=\frac{gsyz}{\theta_0}.
\end{align*}
Here $s<0$ is a constant, so $\bar{\theta}$ represents the decrease in potential temperature when moving away from the equator and towards the north pole.
Observe that $\phi_0+\bar{\phi}$ is in hydrostatic balance with $\theta_0 + \bar{\theta}$, which means that
\begin{align*}
\partial_z(\phi_0+\bar{\phi})=\frac{g}{\theta_0} (\theta_0+\bar{\theta}).
\end{align*}

We seek \emph{vertical slice solutions} of \eqref{eqn:EulBous}-\eqref{eqn:EulBous:3}, where the velocity depends only on $x$, $z$ and $t$. We consider the fluid equations in the $(x,z)$-domain 
\begin{equation}
\label{eqn:domain}
\Omega:=[-L,L)\times [-H/2,H/2],
\end{equation}
we impose $2L$-periodic boundary conditions in $x$, and we require that $w(x,z,t)=0$ for $z\in\{-H/2,H/2\}$.
Substituting 
$\bu=\bu(x,z,t)$, \eqref{eqn:thDe}, and \eqref{eqn:phiDe} into \eqref{eqn:EulBous}-\eqref{eqn:EulBous:3} yields the following form of Eady-Boussinesq vertical slice equations:
\begin{align}\label{eqn:EadySlice1}
	\left\{
	\begin{array}{l@{}c}
		D_t u - fv = -\partial_x\phi,\vspace{3pt}\\
		D_tv + fu = -\frac{gsz}{\theta_0},\vspace{3pt}\\
		D_tw = -\partial_z\phi + \frac{g\theta}{\theta_0},\vspace{3pt}\\
		D_t\theta + s v = 0,\vspace{3pt}\\
		\partial_x u + \partial_z w = 0,
	\end{array}\right.
\end{align}
where $D_t:=\partial_t+u\partial_x+w\partial_z$ is the in-slice material derivative operator.

The \emph{semi-geostrophic approximation} is obtained from \eqref{eqn:EadySlice1} by neglecting the derivatives $D_tu$ and $D_tw$. This results in the system 
\begin{align}\label{eqn:SGES1}
	\left\{
	\begin{array}{l@{}c}
		D_t v + fu = -\frac{gsz}{\theta_0},\vspace{3pt}\\
		D_t\theta + s v = 0,\vspace{3pt}\\
		\partial_x u + \partial_z w = 0,
	\end{array}\right.
\end{align}
where the meridional velocity $v$ and potential temperature $\theta$ are determined by the {geopotential} $\phi$:
\begin{equation}
\label{eqn:vThRel0}
v = \frac{1}{f}\partial_x\phi, \qquad
\theta=\frac{\theta_0}{g}\partial_z\phi.
\end{equation}
In what follows, we refer to this system as the {\em SG Eady slice equations}.

Define the background shear velocity
\begin{align}
\label{eq:ubar}
	\overline{u}(z)=-\frac{gsz}{f\theta_0}.
\end{align}
The
system \eqref{eqn:SGES1}, \eqref{eqn:vThRel0} has steady state 
\begin{align}
\label{eq:SteadyState}
(u,v,w,\theta,\phi)=
\left( \overline{u}(z),0,0, \tfrac{N^2 \theta_0}{g}\left(z+\tfrac{H}{2}\right),
\tfrac{N^2}{2}\left(z+\tfrac{H}{2}\right)^2 \right),
\end{align}
where the constant $N$ is the 
Brunt-V\"{a}is\"{a}l\"{a} or buoyancy frequency. The additive constant in the definition of $\theta$ is chosen so that that the potential temperature is zero on the bottom of the domain, namely $\theta(-H/2)=0$.
We perform a linear instability analysis of this steady state in Appendix \ref{Sec:LinearInstabilityAnalysis} and illustrate its stability numerically in Sections
\ref{Subsec:UNM} and \ref{Subsec:SNM}.

In contrast to the derivations of 
\eqref{eqn:EadySlice1}
given in \cite{visram2014framework} and \cite{yamazaki2017vertical}, in our work $\bar{\theta}$ does not depend on $z$.
This aids the derivation of the SG approximation. 
Instead, the linear dependence on $z$ is built into the steady state for $\theta$.
Also, for notational convenience, the domain $\Omega$ in this paper is a vertical translate of that used in \cite{visram2014framework} and \cite{yamazaki2017vertical}, so $\overline{u}$ differs by an additive constant from the background velocity used therein.

\subsection{Notation}
In what follows, we make a departure from the notation used thus far and in previous works. This is so that the geometric method and its relation to optimal transport can be stated with clarity. We denote by $\bx=(x_1,x_2)$ a position vector in $\Omega$, thereby replacing $(x,z)$ by $\bx$. Likewise, we replace the horizontal and vertical velocities $(u,w)$ by $\bu=(u_1,u_2)$, we let $\nabla=(\partial_{x_1},\partial_{x_2})$ denote the 2-dimensional gradient, and we let $D_t=(\partial_t + \bu\cdot\nabla)$ denote the in-slice material derivative.

\subsection{SG Eady slice model} 
\label{sec:Lagrangian}
Define the modified geopotential $P:\Omega\to \R$ by
\begin{align*}
	P(\mathbf{x},t)=\frac{1}{2}x_1^2+\frac{1}{f^2}\phi(\mathbf{x},t).
\end{align*}
 Using \eqref{eqn:vThRel0} we obtain the relation
\begin{align}\label{eqn:vThP}
	\nabla P(\bx,t)=
\begin{pmatrix}
\frac{1}{f}v(\mathbf{x},t)+x_1
\\
\frac{g}{f^2\theta_0} \theta(\mathbf{x},t)
\end{pmatrix}.
\end{align}
The first component of $\nabla P$ is commonly referred to as the \emph{absolute momentum}.
The semi-geostrophic (SG) Eady slice equations \eqref{eqn:SGES1} then become
\begin{align}\label{eqn:SGESEul}
	\left\{
		\begin{array}{l@{}c}
			\left(\partial_t+\mathbf{u}\cdot\nabla\right)\nabla P= J\big(\text{id}_{\bx}-\left(\nabla P\cdot\mathbf{e}_1\right)\mathbf{e}_1\big),\\
			\nabla\cdot \mathbf{u} = 0,
		\end{array}\right.
\end{align}
where $\text{id}_{\bx}(\bx,t)=\bx$, $\mathbf{e}_1:=(1,0)$, and
\begin{align*}
	J=\frac{gs}{f\theta_0}
	\begin{pmatrix}
	0 & -1 \\
		1 & 0
	\end{pmatrix}.
\end{align*}
The global-in-time existence of solutions of \eqref{eqn:SGESEul} and its 3-dimensional analogue for physically-relevant initial data remains an open problem.
{For future reference, we record that the steady state \eqref{eq:SteadyState} 
corresponds to 
\begin{equation}
\label{eq:Psteady}
\nabla \overline{P}(\mathbf{x})
= 
\begin{pmatrix}
x_1 
\\
\frac{N^2}{f^2} \left( x_2 + \frac H2 \right)
\end{pmatrix}.
\end{equation}
}

We now introduce the Lagrangian form of the equations. Let $\mathbf{F}$ denote the in-slice flow corresponding to $\bu$ so that $\mathbf{F}(\bx,t)\in\Omega$ is the position at time $t$ of a particle that started at position $\bx$, that is
\begin{align*}
\partial_t\mathbf{F}(\mathbf{x},t)=\bu(\mathbf{F}(\mathbf{x},t),t),\qquad \mathbf{F}(\bx,0)=\bx.
\end{align*}
Since $\bu(\cdot,t)$ is divergence free, the flow $\mathbf{F}$ is mass-preserving. Written in terms of the Lagrangian variable
\begin{align}
\label{eq:Z}
	\mathbf{Z}(\mathbf{x}, t):=\nabla P(\mathbf{F}(\mathbf{x}, t), t),
\end{align}
the transport equation \eqref{eqn:SGESEul} becomes
\begin{align}\label{eqn:SGESLagCts}
	\partial_{t} \mathbf{Z}=J\big(\mathbf{F}- \left(\mathbf{Z}\cdot \mathbf{e}_1\right)\mathbf{e}_1\big).
\end{align}
In contrast to the Eulerian setting, global-in-time solutions of the equations in Lagrangian coordinates are known to exist for a wide and physically-relevant class of initial data \cite{feldman2015semi}.

\subsection{Energy}
We define the \emph{total geostrophic energy} at time $t$ to be
\begin{align}\label{eqn:totalGeoEnergy}
	\mathcal{E}(t):=\mathcal{K}_v(t)+\mathcal{P}(t),
\end{align}
where
\begin{align}\label{eqn:KE}
	\mathcal{K}_v(t):=\frac{1}{2}\int_{\Omega}v^2(\bx,t)\,\rd\bx
\end{align}
is the total kinetic energy due to $v$
and
\begin{align}\label{eqn:PE}
	\mathcal{P}(t):=\int_{\Omega}-\frac{g\theta(\bx,t)x_2}{\theta_0}\,\rd\bx
	  +\int_{\Omega}N^2\left(x_2+\frac{H}{2}\right)x_2\,\rd\bx
\end{align}
is the total potential energy.
(The final integral in the definition of $\mathcal{P}(t)$ 
ensures that the total potential energy of the steady state \eqref{eq:SteadyState} is zero. We include it for consistency with \cite{yamazaki2017vertical}.) The total geostrophic energy $\mathcal{E}$ is conserved by the Eulerian equations \eqref{eqn:SGESEul}. 
Note that in contrast to the energy defined in \cite[equations (23)-(26)]{yamazaki2017vertical}, the total geostrophic energy $\mathcal{E}$ does not include the kinetic energy due to the in-slice velocity $\bu$.

\section{The geometric method}\label{sect:geomMethod}
{The geometric method is a spatial discretisation of 
 equation \eqref{eqn:SGESLagCts}, which yields 
 a system of ordinary differential equations. It is derived using an energy minimisation principle known as {\em Cullen's stability principle}.}
It was first described and implemented in \cite{cullen1984extended}, and later used in the context of the SG Eady slice problem in \cite{culrou}.
 {We describe the original implementation in Section \ref{Subsec:Comparison}. In this section we recast the geometric method in the language of \emph{semi-discrete optimal transport theory}, which underpins our novel implementation and will aid its description.}

\subsection{The stability principle and semi-discrete optimal transport}\label{sect:Stability}
In the present context, the stability principle can be stated as follows:
\begin{quotation}
	Stable solutions of \eqref{eqn:SGESLagCts} are those which, at each time, minimise the total geostrophic energy \eqref{eqn:totalGeoEnergy} over all periodic mass-preserving rearrangements of fluid particles that conserve the absolute momentum and potential temperature.
\end{quotation}
This can be shown to be equivalent to assuming that $P$ is convex. 

{In the geometric method, we approximate the modified geopotential $P$ at each time by a piecewise affine function, and apply the stability principle.
{The gradient of} any piecewise affine function on the fluid domain $\Omega$ is uniquely identified by a tessellation of $\Omega$ by cells $S_i$, $i\in\{1,\ldots,n\}$, and corresponding points $\z_i\in\R^2$,  where $\z_i$ is the gradient of the piecewise affine function on $S_i$.}
Each periodic rearrangement of fluid particles corresponds to a different tessellation of $\Omega$ by cells $\widetilde{S}_i$. Such a rearrangement is mass-preserving if $\widetilde{S}_i$ has the same area as $S_i$ for all $i\in\{1,\ldots,n\}$. The absolute momentum and potential temperature are conserved when the corresponding image points $\z_i$ are fixed. 

By writing the geostrophic energy \eqref{eqn:totalGeoEnergy} in terms of the points $\z_i$ and sets $S_i$, accounting for the periodic boundary conditions on $\phi$, and neglecting terms that are constant over all periodic mass-preserving rearrangements, 
the stability principle can be rephrased as a \emph{semi-discrete optimal transport problem} (Definition \ref{def:SDOT}); {see Appendix \ref{app:OTderivation} for the derivation.}
This is a type of optimal partitioning problem.
{We refer to its solution} as an \emph{optimal partition}.

In what follows, for $A\subset \R^2$ we denote by $|A|$ the area of $A$.
 
\begin{defn}[Optimal partition]\label{def:SDOT}
	Given \emph{seeds} $\z=(\z_1,\ldots,\z_n)\in\R^{2n}$ {with $\z_i \ne \z_j$ if $i \ne j$} and \emph{target masses} $\obm=(\om_{1},\ldots,\om_n) { \in \R^n}$ {with $\om_{i}>0$ and $\sum_{i=1}^n\om_{i}=|\Omega|$}, a partition of $\Omega$ is said to be optimal if it minimises the \emph{transport cost}
	\begin{align}\label{eqn:transportCost}
		\mathcal{T}(\{S_i\}_{i=1}^n):=\sum_{i=1}^n\int_{S_i}|\bx- \z_i|_{\per}^2\, \rd \bx
	\end{align}
	{among} all partitions $\{S_i\}_{i=1}^n$ of $\Omega$ that satisfy the mass constraint
	\begin{align}\label{eqn:massConstraint}
		|S_i|=\om_{i}\quad \forall\, i\in\{1,\ldots,n\}.
	\end{align}
	Here
	$|\bx-\by|_{\per}$ is the distance between $\bx, \by \in\R^2$ taking into account $2L$-periodicity in the first component:
	\begin{align*}
		|\bx-\by|_{\per}:=\min_{\mathbf{k}\in K}|\bx-\by-\mathbf{k}|,
	\end{align*}
where 
	\begin{align*}
		K:=\left\{2Lk\mathbf{e}_1\,\vert\, k\in\Z\right\}.
	\end{align*}
\end{defn}

\revB{Recall that the characteristic function $\mathds{1}_{A}:\R^2\to\R$ of a set $A\subset\R^2$ is defined by
\begin{align*}
\mathds{1}_{A}(\bx):=\left\{
\begin{array}{l}
\, 1 \quad \bx \in A,\\
\, 0 \quad \text{otherwise.}
\end{array}
\right.
\end{align*}}Any modified geopotential {$P$} that {is piecewise affine} at time $t$ and satisfies the stability principle has the form
\begin{align}\label{eqn:pwConstNablaP}
	\nabla P (\bx,t) = \sum_{i=1}^n \left(\z_i	+\mathbf{k}_*(\bx,\z_i)\right)\mathds{1}_{S_i}(\bx),
\end{align}
where {$\{S_i(t)\}_{i=1}^n$} is an optimal partition corresponding to the seeds {$\z_i(t)$} in the sense of Definition \ref{def:SDOT}, and
\begin{align}\label{eqn:kTranslate}
\mathbf{k}_*(\bx,\z_i):=\underset{\mathbf{k}\in K}{\text{argmin}}\vert\bx-\z_i -\mathbf{k}\vert
\end{align}
 accounts for the periodic boundary condition on the geopotential $\phi$ (cf.~\cite[Theorem 1.25]{santambrogio2015optimal}). 
{Note that $\mathbf{k}_*(\bx,\z_i)$ is well defined for almost-every $\bx\in\Omega$.}
Moreover, without loss of generality, ${\z_i(t)} \in [-L,L)\times \R$ for each $i\in\{1,\ldots,n\}$.

Optimal partitions can be described in terms of periodic Laguerre diagrams. These are partitions of the domain $\Omega$ into cells parametrised by the seeds $\z=(\z_1,\ldots,\z_n)$ and a set of weights $\w=(w_1,\ldots,w_n)\in\R^n$. 

\begin{defn}[Periodic Laguerre diagram]\label{def:LagCells}
Let $\z=(\z_1,\ldots,\z_n)\in\left([-L,L)\times\R\right)^{n}$ {with $\z_i \ne \z_j$ if $i\ne j$. Let} $\w=(w_1,\ldots,w_n)\in\R^n$. \revB{For $i\in\{1,\ldots,n\}$, we define the set}
\begin{align}\label{eqn:LagIneq}
	C_{i,\per}(\z,\w)
	:=
	\left\{
	\bx\in \Omega  :
	{|\bx-\z_i|_{\per}^2-w_i\leq |\bx-\z_j|_{\per}^2-w_j
	\; \forall \, j\in\{1,\ldots,n\}}
	\right\}.
\end{align}
\revB{This is the $i$\textsuperscript{th} periodic Laguerre cell generated by $(\z,\w)$, and the collection of all cells $\{C_{i,\per}\}_{i=1}^n$ is the periodic Laguerre diagram generated by $(\z,\w)$.}
\end{defn}

\begin{defn}[Cell-area map]\label{defn:cellAreas}
	We define the cell-area map  $\m=(m_1,\ldots,m_n):\R^{2n}\times \R^n\to\R^n$ by
	\begin{align*}
		m_i(\z,\w):=\left\vert C_{i,\per}(\z,\w)\right\vert.
	\end{align*}
\end{defn}
It is well known that {given} $\z=(\z_1,\ldots,\z_n)\in\R^{2n}$, {with} $\z_i\neq\z_j$ whenever $i\neq j$, {and given
$\obm=(\om_{1},\ldots,\om_n) \in \R^n$ with $\om_{i}>0$,  $\sum_{i=1}^n\om_{i}=|\Omega|$},
there exists a unique weight vector $\w_*(\z)\in\R^n$ with final entry $0$ that satisfies the mass constraint
\begin{align*}
	\m(\z,\w_*(\z))=\obm
\end{align*}
(see for example \cite[Corollary 39]{merigot2021optimal}). The optimal partition of $\Omega$ is then the periodic Laguerre diagram generated by $(\z,\w_*(\z))$. We call $\w_*(\z)$ the {\em optimal weight vector }for the seeds $\z$. Observe that, for all $\lambda\in\R$,
\begin{align}
	C_{i,\per}(\z,\w+\lambda\mathbf{e})=C_{i,\per}(\z,\w),
\end{align}
where $\mathbf{e}=(1,\ldots,1)\in\R^n$. Choosing the last entry of $\w_*(\z)$ to be $0$ ensures uniqueness of the optimal weight vector, without loss of generality. 

To solve the semi-discrete optimal transport problem (Definition \ref{def:SDOT}), it is therefore sufficient to compute the optimal weight vector. 
To do this, we use the well-known fact that $\w_*(\z)\in\R^n$ is the unique maximum (with final entry $0$) of the {\em Kantorovich functional}, 
which is the concave function $\mathcal{K}:\R^n\to\R$  defined by
\begin{align}
	\label{eqn:KantFunct}
	\mathcal{K}(\w):=\sum_{i=1}^n{\int_{C_{i,\per}(\z,\w)}} |\bx-\z_i|_{\per}^2\,\rd\mathbf{x}
	+ \w\cdot (\obm-\m(\z,\w)).
\end{align}
{See for example \cite[Theorem 40]{merigot2021optimal}}.
In effect, the constrained optimal {partitioning} problem is transformed into an unconstrained, finite{-}dimensional, {concave} maximisation problem, which is numerically tractable.

\subsection{Spatial discretisation}\label{sect:ODEderivation}
In summary, we seek solutions of \eqref{eqn:SGESLagCts} for which the associated geopotential $P$ is piecewise affine in space at each time $t$. By the stability principle, $\nabla P$ must have the form \eqref{eqn:pwConstNablaP} with
\begin{align}\label{eqn:ansatz}
	S_i=C_{i,\per}(\z,\w_*(\z))
\end{align}
for some time-dependent map $\z=(\z_1,\ldots,\z_n)$,
where
$\w_*$
is the maximum of $\mathcal{K}$ {and the}
chosen target masses $\obm=(\om_{1},\ldots,\om_n)$ do not change over time. The assumption that $S_i$ has the form \eqref{eqn:ansatz} is equivalent to assuming that $P$ is convex.


\revA{We now sketch the derivation of an ODE for $\z$, further details of which can be found in Appendix \ref{app:ODEderivation}. We begin by making the following definition, which arises naturally in the derivation due to the periodic boundary conditions.}

\begin{defn}
\label{def:Laguerre}
	Let $\z=(\z_1,\ldots,\z_n)\in\R^{2n}$ {with $\z_i \ne \z_j$ if $i\ne j$}. {Let} $\w=(w_1,\ldots,w_n)\in\R^n$. 
	For $i\in\{1,\ldots,n\}$, {we define the convex polygon $C_i(\z,\w)$ by}
	\begin{align*}
		 C_{i}(\z,\w):=
		 \left\{\bx\in\R\times\left[-\frac H2,\frac H2\right] : 
		 {|\mathbf{x}-\z_i|^2-w_i\leq |\mathbf{x}-\z_j-\mathbf{k}|^2-w_j
		 \;
		 \forall \, j \in\{1,\ldots,n\}, 
		 \,
		 k\in K}
		 \right\}.
	\end{align*}
\end{defn}

{Note the difference between Definitions \ref{def:LagCells} and \ref{def:Laguerre}: Definition \ref{def:LagCells} defines a \emph{periodic} Laguerre tessellation of $[-L,L)\times [-H/2,H/2]$ generated by a finite set of seeds and weights $(\z,\w)$, while Definition \ref{def:Laguerre} defines a \emph{non-periodic} Laguerre tessellation of $\mathbb{R} \times [-H/2,H/2]$ generated by all periodic copies of $(\z,\w)$, namely by
$(\z_i+\mathbf{k},w_i)$ for all $i \in \{1,\ldots,n\}$, $\mathbf{k} \in K$.}
In general, the cells $C_i(\z,\w)$ are not contained in $\Omega$, but if $C_i(\z,\w)\subset\Omega$ then $C_i(\z,\w)=C_{i,\per}(\z,\w)$. In any case, both $C_i(\z,\w)$ and $C_{i,\per}(\z,\w)$ have the same area for each $i\in\{1,\ldots,n\}$.

By definition of the optimal weight vector, at each time $t$ the mass constraint
\begin{align*}
	\m\Big(\z(t),\w_*\big(\z(t)\big)\Big)=\obm
\end{align*}
is satisfied. For brevity, we define
\begin{equation}
    \label{eq:brevity}
C_{i,\per}(t)=C_{i,\per}(\z(t),\w_*(\z(t))).
\end{equation}
{Assume} that there exists a mass-preserving flow $\mathbf{F}$ such that
\begin{align}\label{eqn:cellFlow}
	C_{i,\per}(t)=\mathbf{F}\left(C_{i,\per}(0),t\right)
\end{align}
at all times $t$. For each $i\in\{1,\ldots,n\}$, we define $\bc_i(\z)$ to be the centroid of $C_{i}(\z,\w_*(\z))$, namely
\begin{align}
\label{eq:c}
	\bc_i(\z):=\frac{1}{|C_{i}(\z,\w_*(\z))|}\int_{C_{i}(\z,\w_*(\z))}x\ \rd x,
\end{align}
and we define
\begin{align*}
	\bc(\z):=(\bc_1(\z),\ldots,\bc_n(\z)).
\end{align*}
Formally, after substituting the ansatz {\eqref{eqn:pwConstNablaP}}, \eqref{eqn:ansatz} into \eqref{eqn:SGESLagCts}, integrating over the cell {$C_{i,\per}(0)$},
and using \eqref{eqn:cellFlow} to {apply} the change of variables $\mathbf{x}\mapsto \mathbf{F}(\mathbf{x},t)$, one finds that the vector of seed trajectories $\z$ {satisfies} the ODE
\begin{align}
\label{eq:ODE1}
    \left\{
	\begin{array}{l@{}c}
	\dot{\z}_i=J\big(\bc_i(\z)-\left(\z_i\cdot \mathbf{e}_1\right)\mathbf{e}_1\big),\\
	\z_i(0)=\obz_{i}.
	\end{array}
	\right.
\end{align}
for all $i\in\{1,\ldots,n\}$,
{where $\obz=(\obz_1,\ldots,\obz_n)$ denotes the initial seed positions: see Appendix \ref{app:ODEderivation} for the details of the derivation.} {The ambient space containing individual seed positions $\z_i(t)$ is referred to as \emph{geostrophic space}.}
The ODE \eqref{eq:ODE1} is the discrete analogue of the Lagrangian equation \eqref{eqn:SGESLagCts}, where the centroid map $\bc$ plays the role of the flow $\mathbf{F}$. 
The full system can be written compactly as
\begin{align}\label{eqn:ODE}
    \left\{
	\begin{array}{l@{}c}
	\dot{\z}=\mathbf{J}\left(\bc(\z)-\mathbf{\Pi}\z\right),\\
	\z(0)=\obz,
	\end{array}
	\right.
\end{align}
where $\mathbf{J}$ and $\mathbf{\Pi}$ are the $2n\times 2n$ matrices
\begin{align*}
	\mathbf{J}=\text{diag}(J,\ldots,J) \quad \text{and} \quad
	\mathbf{\Pi}=\text{diag}(\mathbf{e}_1\otimes \mathbf{e}_1,...,\mathbf{e}_1\otimes \mathbf{e}_1).
\end{align*}
The matrix $\mathbf{\Pi}$ acts on each seed $\z_i$ by projection onto the horizontal coordinate. We will use this form of the system when we describe the time discretisation of the system: see Algorithm \ref{alg:adaptiveTimeStep} below. 

\subsection{Structure preservation and recovery of physical variables}\label{sect:Precovery}

If $\z$ satisfies the ODE \eqref{eqn:ODE} and 
{$\w(t)=\w_*(\z(t))$,}
then the total geostrophic energy \eqref{eqn:totalGeoEnergy} corresponding to a modified geopotential $P$ of the form \eqref{eqn:Precovery} is constant in time. The geometric method therefore inherits the energy conservation property possessed by the Eulerian SG Eady slice equation \eqref{eqn:SGESEul}. 
{Moreover, the one-parameter family of Laguerre tessellations $t\mapsto \{C_i(\z(t),\w(t))\}_{i\in\{1,\ldots,n\}}$
corresponds to an area-preserving flow of the fluid since 
the mass of each cell is conserved by definition of $\w_*$.}

Given seed trajectories $\z=(\z_1,\ldots,\z_n)$ satisfying the ODE \eqref{eqn:ODE}, let
{$\w(t)=\w_*(\z(t))$.}
The corresponding modified geopotential $P$ is
\begin{align}
    \label{eqn:Precovery}
	{P(\mathbf{x},t)=}
	\sum_{i=1}^n&\Big(\big(\z_i+{\mathbf{k}_*(\bx,\z_i)}\big)\cdot\mathbf{x}
	-\frac{1}{2}\left\vert\z_i+{\mathbf{k}_*(\bx,\z_i)}\right\vert^2
	+\frac{1}{2}w_i\Big)\mathds{1}_{C_{i,\per}(\z,\w)}(\bx),
\end{align}
{and $\mathbf{k}_*(\bx,\z_i)$ was defined in \eqref{eqn:kTranslate}.} \revA{This is a piecewise affine convex function whose gradient is defined almost-everywhere and is given by}
\begin{align}
\label{eq:nablaP}
	\nabla P(\mathbf{x},{t}) = \sum_{i=1}^n\big(\z_i+{\mathbf{k}_*(\bx,\z_i)}\big)\mathds{1}_{C_{i,\per}(\z,\w)}(\bx).
\end{align}
Expressions for the approximate out-of-slice velocity {and} potential temperature can then be recovered using \eqref{eqn:vThP}.

\section{Implementation}\label{sect:implementation}

We now give details of our numerical implementation of the geometric method. \revC{Our software is available at the following GitHub repositories.}
\vspace{0.2cm}

\noindent \revC{\emph{SG-Eady-Slice}: MATLAB functions for initialising and solving the SG Eady slice equations with periodic boundary conditions in the horizontal direction using the geometric method of Cullen \& Purser \cite{cullen1984extended} with our adaptive time-stepping scheme.\newline
\texttt{\url{https://github.com/CharliePEgan/SG-Eady-slice}}}
\vspace{0.2cm}

\noindent \revC{\emph{MATLAB-Voro}: MATLAB mex files for generating $2$D and $3$D periodic and non-periodic Laguerre tessellations using Voro++ \cite{Voro++}.\newline
\texttt{\url{https://github.com/smr29git/MATLAB-Voro}}}

\subsection{Generating discrete initial data}
\label{subsec:quantizeIC}
For fixed $n\in \N$, we approximate given initial data
{$\nabla P(\mathbf{x}, 0)$}
by a piecewise constant function $\nabla P^n_0$ of the form \eqref{eq:nablaP}
defined by initial seeds $\obz=(\obz_1,\ldots,\obz_N)$ and cell areas $\obm=(\om_1,\ldots,\om_n)$. {We now describe how we choose $\obz$ and $\obm$.} 

{In Section \ref{sect:results}
we will take $\nabla P(\mathbf{x}, 0)$ to be a small perturbation of the steady state $\nabla \overline{P}$, which was defined in \eqref{eq:Psteady}. Therefore}
 the seeds $\obz_i$ are taken from the image domain $\nabla P(\Omega,0)$, which
is a small perturbation of the rectangle 
\begin{align*}
	R=[-L,L)\times\left[0,\frac{N^2H}{f^2}\right].
\end{align*}
{We define $\obz_i$ as follows:}
\begin{enumerate}[leftmargin=*]
	\item use Lloyd's algorithm \cite{DuFaberGunzburger1999} to generate points $\oby=(\oby_1,\ldots,\oby_n)$ that are approximately uniformly distributed in $R$;
	\item for each $i\in\{1,\ldots,n\}$, map $\oby_i$ into $\Omega$ by 
	{inverting $\nabla \overline{P}$}:
	\[
	\obx_i:=
	{\left(\nabla \overline{P}\right)^{-1}(\oby_i)} =
	\begin{pmatrix}
		1 & 0 \\
		0 & \frac{f^2}{N^2}
	\end{pmatrix}
	\oby_i
	{- 
	\begin{pmatrix}
	0 \\ \frac H2
	\end{pmatrix}};
	\]
	\item define $$\obz_i:=\nabla P(\obx_i,0).$$
\end{enumerate}

We briefly describe Lloyd's algorithm, which is an iterative method for quantising measures (see, for example, \cite[Section 5.2]{DuFaberGunzburger1999}). Let $\by=(\by_1,\ldots,\by_n)\in\left([-L,L)\times\R\right)^{n}$ with $\by_i \ne \by_j$ if $i\ne j$. For each $i\in\{1,\ldots,n\}$, \revB{define the $i$\textsuperscript{th} periodic Voronoi cell generated by $\by$ as}
\begin{align}
    V_{i,\per}(\by):=
    \big\{
    \mathbf{r} \in R : 
    {|\mathbf{r}-\by_i|_{\per}\leq |\mathbf{r}-\by_j|_{\per}}
    \; \forall j\in \{1,\ldots, n\}
    \big\}.
\end{align}
\revB{The collection of all cells $\{V_{i,\per}(\by)\}_{i=1}^n$ is called the periodic Voronoi tessellation of $R$.}
At each iteration of the Lloyd algorithm, the value of $\by_i$ is updated by moving it to the centroid of its Voronoi cell:
\begin{align*}
	\by_i\mapsto \frac{1}{\left\vert V_{i,\per}(\by)\right\vert}\int_{V_{i,\per}(\by)}\big(\mathbf{r}-\mathbf{k}_*(\mathbf{r},\by_i)\big)\,\rd\mathbf{r}
\end{align*}
for each $i\in\{1,\ldots,n\}$.
The centroids can be computed exactly (without numerical integration): see for example \cite[online supplementary material, equation (4)]{bourne2015centroidal}.
In our implementation of Lloyd's algorithm, we start with points $\by=(\by_1,\ldots,\by_n)$ on a regular triangular lattice in $R$ and perform $100$ iterations to obtain $\oby=(\oby_1,\ldots,\oby_n)$. The target areas are defined for each $i\in\{1,\ldots,n\}$ as
\begin{align*}
	\om_i=\frac{f^2}{N^2}\left\vert V_{i,\per}(\oby)\right\vert.
\end{align*}

The method above generates points $\obz_i$ that are approximately optimally sampled from the distribution 
on $\nabla P(\Omega,0)$ that is obtained by pushing forward the uniform distribution on $\Omega$ by $\nabla P(\cdot,0)$.
It is numerically cheaper than directly sampling this distribution using Lloyd's algorithm since it avoids numerical integration.

{Since the discrete dynamics \eqref{eqn:ODE} conserves the total geostrophic energy $\mathcal{E}$, it is desirable that the initial condition starts on the correct energy surface, i.e. that the discrete initial data $\nabla P_0^n$ has the same total geostrophic energy as $\nabla P(\mathbf{x},0)$. While there exists discrete initial data with this property \cite{egan2022constrained}, the choice above does not. Nevertheless, it is easy to generate and it approximates the initial energy well enough if $n$ is sufficiently large.}

{
\subsection{Generating Laguerre tessellations}
\label{subsec:SolveOT}
In two dimensions the worst-case complexity of computing a Laguerre tessellation with $n$ seeds is $\mathcal{O}(n \log n)$, where $n$ is the number of seeds. This can be achieved for example by the lifting method of Aurenhammer \cite{Aurenhammer1987}. We computed Laguerre and Voronoi tessellations using the C++ library Voro++ \cite{Voro++} and our own mex file to interface with MATLAB. {While this does not achieve the optimal scaling in $n$, it is sufficiently fast for the values of $n$ that we use.} 
}

\subsection{Solving the semi-discrete optimal transport problem}

In any numerical scheme for solving the ODE \eqref{eqn:ODE}, it is necessary to evaluate the right-hand side at each 
{time step}. Given seeds $\z$, this involves solving the corresponding semi-discrete optimal transport problem to compute $\bc(\z)$.
{This is the most expensive part of the algorithm. Here we do this by finding}
the maximum of the Kantorovich functional $\mathcal{K}$ (equation \eqref{eqn:KantFunct}) using the damped Newton method  developed in \cite{kitagawa2016convergence}; {see Algorithm \ref{alg:DN}. This solves the nonlinear algebraic equation $\nabla \mathcal{K} = 0$.}

First, a guess $\w$ for the optimal weight vector is proposed. The Newton direction $\mathbf{d}$ is then determined by solving the linear system
\begin{align}\label{eqn:NewtDir}
	\left\{\begin{array}{l}
		\mathrm{D}^2\mathcal{K}(\w)\mathbf{d}={-}\nabla \mathcal{K}(\w), \\
		\mathbf{d}\cdot \mathbf{e}_n=0,
	\end{array}\right.
\end{align}
where $\mathbf{e}_n=(0,\ldots,0,1)\in\R^n$. (Formulas for $\mathrm{D}^2\mathcal{K}$ and $\nabla \mathcal{K}$ are given in Appendix \ref{app}.) For the linear system \eqref{eqn:NewtDir} to have a unique solution, it is necessary and sufficient for $\w$ to satisfy the mass-positivity condition
\begin{align}\label{eqn:massPositivity}
	\left\vert C_{i,\per}(\z,\w)\right\vert>0,
\end{align}
for all $i\in\{1,\ldots,n\}$. To ensure that this condition is met by the subsequent iterate, backtracking is used to determine the length of the Newton step. The algorithm terminates once 
{$|\nabla \mathcal{K}|$ is less than a given tolerance.}

Convergence of the damped Newton algorithm is guaranteed if and only if the initial guess for the weights satisfies \eqref{eqn:massPositivity}, in which case it converges globally with linear speed, and locally with quadratic speed, as the number of iterations diverges (see \cite[Proposition 6.1]{kitagawa2016convergence}). Our adaptive time-stepping scheme (Algorithm \ref{alg:adaptiveTimeStep}) for solving the ODE \eqref{eqn:ODE} provides a robust way to generate a good {initial} guess for the weights at each 
{time step}
given the optimal weights from the previous 
{time step}.

It remains to generate a good {initial} guess for the weights at time $t=0$. 
By \cite[Theorems 3 and 4]{levy2018notions}, each cell $C_{i,\per}(\z,\w)$ is non-empty if and only if \revC{$\w$ is $c$-concave in that sense that }there exists \revC{$\varphi:\Omega\to\R$} such that
\begin{align}\label{eqn:ctrans}
	w_i=\min_{\bx\in\Omega}\{|\bx-\z_i|_{\per}^2-\varphi(\revC{\bx})\}.
\end{align}
{This condition ensures that the cells $C_{i,\per}(\z,\w)$ are non-empty but not necessarily that they satisfy the mass-positivity condition \eqref{eqn:massPositivity}.}
{If $\varphi=0$, then $\w$ defined by \eqref{eqn:ctrans} satisfies \eqref{eqn:massPositivity} 
if in addition}
the horizontal components of the seeds $\z_i$ are distinct {(or if $\z_i \in \Omega$ for all $i$, but this is not the case for our simulations in Section \ref{sect:results})}. 
{At time $t=0$, we therefore first randomly perturb the horizontal components of the seeds:}
\[
\widetilde{\z}_i:=\obz_i+\zeta_i {\mathbf{e}_1},
\]
where $\zeta_i\in\R$ are appropriately scaled, randomly generated numbers.
{Then we}
apply Algorithm \ref{alg:DN} with seeds $\widetilde{\z}_i$ and initial weights
\begin{align}\label{eqn:ctrans0}
	w_i=\min_{\bx\in\Omega}|\bx-\widetilde{\z}_i|_{\per}^2
\end{align}
to obtain $\w_*(\widetilde{\z})$. Since $\w_*$ is continuous, for sufficiently small perturbations $\zeta_i$, $\w_*(\widetilde{\z})$ is a good {initial} guess for $\w_*(\obz)$ (which is computed in the initial step of Algorithm \ref{alg:adaptiveTimeStep}).
Note that \eqref{eqn:ctrans0} is the square of the $K$-periodic distance from $\z_i$ to $\Omega$ which can be simply computed as
\begin{align*}
	{w_i =}
	\begin{cases}
		\left(\widetilde{\z}_i\cdot\mathbf{e}_2-\frac{H}{2}\right)^2 \quad &\text{if} \quad \widetilde{\z}_{i}\cdot\mathbf{e}_2>\frac{H}{2},\\
		0 \quad &\text{if}\quad
		{\widetilde{\z}_{i}} \in \Omega, \\
		\left(\widetilde{\z}_{i}\cdot\mathbf{e}_2+\frac{H}{2}\right)^2 \quad &\text{if} \quad \widetilde{\z}_{i}\cdot\mathbf{e}_2<-\frac{H}{2},
	\end{cases}
\end{align*}
where $\mathbf{e}_2=(0,1)$.

\begin{algorithm}
	\caption{Damped Newton algorithm of Kitagawa, M\'{e}rigot and Thibert \cite{kitagawa2016convergence}}
	\label{alg:DN}
	\textbf{Input:} Seeds $\z  =(\z_1,\ldots,\z_n)$, target masses $\obm=(\om_{1},\ldots,\om_{n})$, a percentage mass tolerance $\eta$, and an initial guess for the weights $\w$ such that
	\[
	\varepsilon:=\frac{1}{2} \min \left[\min_{i} m_i(\z,\w), \min_{i} \om_{i}\right]>0.
	\]
	\textbf{Initialisation:} Set $\w^{0}=\w$ and $k=0$. Convert $\eta$ to an absolute mass tolerance:
	\[
	\eta_{\mathrm{abs}}:=\frac{\eta}{100}\min_i{\overline{m}_i}.
	\]
	\\
	\textbf{While:} $\left\|\m\left(\z,\w^{(k)}\right)-\obm\right\|_{\infty} \geqslant \eta_{\mathrm{abs}}$:\\
	\textbf{Step 1:} Solve the following linear system for the Newton direction $\mathbf{d}^{(k)}$:
	\begin{align}\label{eqn:NewtDirAlg}
		\left\{\begin{array}{l}
			\mathrm{D}^2\mathcal{K}\left(\w^{(k)}\right)\mathbf{d}^{(k)}= - \nabla \mathcal{K}\left(\w^{(k)}\right), \\
			\mathbf{d}^{(k)}\cdot \mathbf{e}_n=0.
		\end{array}\right.
	\end{align}
	\textbf{Step 2:} Determine the length of the Newton step using backtracking: find the minimum $\ell \in \mathbb{N}\cup\{0\}$ such that $\w^{(k, \ell)}:=\w^{(k)}+2^{-\ell} \mathbf{d}^{(k)}$ satisfies
	\[
	\left\{
	\begin{array}{l}
		m_i\left(\z,\w^{(k, \ell)}\right) \geqslant \varepsilon \quad \forall\, i\in\{1,\ldots,n\}, \phantom{\Big|}\\
		\left\|\m\left(\z,\w^{(k, \ell)}\right)-\obm\right\|_{\infty} \phantom{\Big|}\leq\left(1-2^{-(\ell+1)}\right)\left\|\m\left(\z,\w^{(k)}\right)-\obm\right\|_{\infty}.
	\end{array}
	\right.
	\]
	\textbf{Step 3:} Define the damped Newton update $\w^{(k+1)}:=\w^{(k)}+2^{-\ell} \mathbf{d}^{(k)}$ and $k \leftarrow k+1$.\\
	\textbf{Output:} A vector $\w^{(k)}$ such that 
	\[
	\frac{100}{\underset{i}{\min}\ \overline{m}_i}\left\|
	\m\left(\z,\w^{(k)}\right)-\obm
	\right\|_{\infty} \leq \eta.
	\]
\end{algorithm}

\subsection{Solving the ODE using an adaptive time-stepping method}
\label{subsec:adaptive}

The most computationally intensive part of solving the ODE \eqref{eqn:ODE} numerically is evaluating the function $\bc$.
In addition, $\bc$ is continuously differentiable (on the set of distinct seed positions), but not twice continuously differentiable in general: see \cite{bourne2022semi}. As such, we can only expect convergence of the ODE solver up to second order in the time step size.
We therefore use an Adams-Bashforth 2-step method (AB2) with adaptive time-stepping. This requires only one function evaluation at each time step.

Our adaptive time-stepping scheme (Algorithm \ref{alg:adaptiveTimeStep}) speeds up the function evaluations by using the ODE \eqref{eqn:ODE} to generate a good initial guess for the weights for Algorithm \ref{alg:DN}, as we now describe.
Suppose that at a given time step we have seeds $\z$ and optimal weights $\w$. Following the notation of Algorithm \ref{alg:adaptiveTimeStep}, for a proposed {time} step size $h^l$, the AB2 scheme determines an increment $\z^{l}_{\text{inc}}$ for the seed positions. The seeds at the subsequent time step are then defined as
\begin{align*}
	\z^{l} := \z + \z^{l}_{\text{inc}}.
\end{align*}
We generate an initial guess $\w^{l}$ for the weights corresponding to the seeds $\z^{l}$ by taking a first order Taylor expansion of $\w_*(\z^l)$ around $\z$:
\begin{align*}
	\w^{l} := \w+D_{\z}\w_*(\z) \, \z^{l}_{\text{inc}}.
\end{align*}
If the seeds and weights $(\z^l,\w^l)$ generate a periodic Laguerre tessellation of $\Omega$ with no zero-area cells, then the proposed step size is accepted. Otherwise, the seeds and weights $(\z^l,\w^l)$ are discarded, the proposed step size is halved and the updates are recalculated. This ensures the convergence of the damped Newton algorithm when used at the subsequent time step.

\begin{algorithm}
	\caption{Adaptive time-stepping scheme}
	\label{alg:adaptiveTimeStep}
	\textbf{Input:} Initial seeds $\obz$ and masses $\obm$, a final time $T$, a default time-step size $h_{\text{def}}$ (seconds), and a percentage mass tolerance $\eta$.\\
	\textbf{Initial step:} Compute the optimal weight vector $\w_*(\obz)$ using Algorithm \ref{alg:DN}, and do a forward Euler step to determine seed positions at time $h_{\text{def}}$. 
	Set
	\setlength{\abovedisplayskip}{4pt}
	\setlength{\belowdisplayskip}{4pt}
		\begin{gather*}
		    h=h_{\text{def}},  \quad \w=\w_*(\obz), \quad t=h,
		    \\
		    \bv_{\text{curr}}  = \mathbf{J}\left(\bc(\obz)-\mathbf{\Pi}\obz\right), \quad 	\z = \obz+h\bv_{\text{curr}}.
		\end{gather*}
	\textbf{While:} $t\leq T$:\\
	\textbf{Step 1:} Compute the optimal weight vector $\w_*(\z)$ using Algorithm \ref{alg:DN} with $\w$ as the initial guess for the weights. Set
		\begin{align*}
			&\w \leftarrow \w_*(\z),\\
			&\left(\bv_{\text{prev}},\bv_{\text{curr}}\right) \leftarrow \left(\bv_{\text{curr}}, \mathbf{J}\left(\bc(\z)-\mathbf{\Pi}\z\right)\right).
		\end{align*}
	\textbf{Step 2:} Determine the minimum $l\in\N\cup\{0\}$ such that $\min_{i\in\{1,...,N\}}(m_i(\z^l,\w^l))>0$, where the updated seeds $\z^l$ and weights $\w^l$ are defined for $l\in\N$ as follows:
	\newline
{
	Propose step size:
	\setlength{\abovedisplayskip}{4pt}
	\setlength{\belowdisplayskip}{4pt}
	\begin{align*}
		h_l \leftarrow \frac{h_{\text{def}}}{2^l}
	\end{align*}
	Compute AB2 coefficients:
    \[
		\left(c_{\text{prev}}^l,c_{\text{curr}}^l\right) \leftarrow \left(-\frac{h_l^2}{2h}, \frac{1}{2}\left(\frac{(h_l+h)^2}{h}-h\right)\right)
	\]
	Compute increments of seeds and weghts:
\begin{gather*}
    		\z_{\text{inc}}^l  \leftarrow c_{\text{prev}}^l\bv_{\text{prev}}+c_{\text{curr}}^l\bv_{\text{curr}}
\\
		\w_{\text{inc}}^l  \leftarrow D_{\z}\w_*(\z)\z_{\text{inc}}^l
			\end{gather*}
	Update seeds and guess for weights:
\[
		\left(\z^l,\w^l\right) \leftarrow \left(\z + \z_{\text{inc}}^l, \w_*(\z) + \w_{\text{inc}}^l\right)
	\]
	}\textbf{Step 3:} Set 
		\begin{align*}
			(\z,\w)\leftarrow (\z^l,\w^l),\quad t\leftarrow t+h_l,\quad h\leftarrow h_l.
		\end{align*}
	\textbf{Output:} Seed positions at each time step.
\end{algorithm}

From numerical experiments, we observed that this method is up to forty times faster than either of the following methods for generating the initial guess for the weights: (i) using the weights from the previous time step, and a small enough time step to ensure that the cells have positive area; (ii) using \eqref{eqn:ctrans0} with a suitable perturbation of the seeds, without using any information from the previous time step or adaptive time stepping.
Moreover, the sensitivity analysis in Section 
\ref{sect:sensitivity}
shows that the adaptively chosen times step sizes are suitable, i.e. not unreasonably small.

The expression for $D_{\z}\w_*(\z)$ can be obtained by implicit differentiation of the mass constraint.
Indeed, by definition of the optimal weight map $\w_*$,
\begin{align}\label{eqn:massConstraintOpt}
	\m(\z,\w_*(\z))=
	\obm.
\end{align}
By differentiating \eqref{eqn:massConstraintOpt} with respect to $\z$, we see that the derivative $D_{\z}\w_*(\z)$ satisfies the equation
\begin{align}\label{eqn:linSystemForDzw}
	D_{\w}\m(\z,\w_*(\z))D_{\z}\w_*(\z)= -D_{\z}\m(\z,\w_*(\z)).
\end{align}
The matrix $D_{\w}\m(\z,\w_*(\z))$ is symmetric and singular with kernel spanned by $\mathbf{e}:=(1,\ldots,1)\in\R^n$.
Each column of $D_{\z}\m(\z,\w_*(\z))$ is orthogonal to $\mathbf{e}$.
Hence, the linear system \eqref{eqn:linSystemForDzw} has a solution $D_{\z}\w_*(\z)$. We choose the solution
\begin{align}\label{eqn:Dzw}
	D_{\z}\w_*(\z)=\left(\begin{array}{c}
		-\mathbf{A}^{-1}\mathbf{B}\\
		\mathbf{0}
		\end{array}
		\right),
\end{align}
where $\mathbf{0} \in \mathbb{R}^{2n}$ is the zero vector and
the matrices $\mathbf{A}\in \R^{(n-1)\times(n-1)}$ and $\mathbf{B}\in\R^{(n-1)\times 2n}$ are defined by
\begin{align}
	\label{eqn:matrixA}
		\mathbf{A}_{i,j}:=\pdone{m_i}{w_j}(\z,\w_*(\z)),
\end{align}
for $i,\,j\in\{1,\ldots,n-1\}$ and
\begin{align}
	\label{eqn:matrixB}
		\mathbf{B}_{k,2l-1}:=\pdone{m_k}{\z_{l,1}},\quad \mathbf{B}_{k,2l}:=\pdone{m_k}{\z_{l,2}},
\end{align}
for $k\in\{1,\ldots,n-1\}$ and $l\in\{1,\ldots,n\}$. Expressions for the derivatives of $\m$ with respect to $\z$ and $\w$ in the non-periodic setting are given for example in \cite{de2019differentiation} and \cite[Lemma 2.4]{bourne2015centroidal}. In Appendix \ref{app}, we state the analogous expressions in the periodic setting. 
	
	\subsection{Comparison with the original implementation}
	\label{Subsec:Comparison}
\revA{In the original presentation of the geometric method \cite{cullen1984extended}, 
a convex piecewise affine modified geopotential $P$ is constructed directly at each time $t$.
The projection of its graph onto the $(x_1,x_2)$-plane
gives a tessellation $\{C_i\}_{i=1}^n$ of $\Omega$. While the term `Laguerre diagram' is not used in \cite{cullen1984extended}, this} tessellation is equivalent to the optimal partition described above in terms of Laguerre cells (up to the inclusion of periodic boundary conditions). The existence and uniqueness of such a map $P$ was also proved formally.

As described in \cite[Section 5.1.2]{cullen2021mathematics}, for some unknown scalars $p_j,\, j\in\{1,\ldots,n\}$, each cell $C_i$ consists of all $\bx\in\Omega$ such that
\begin{align}\label{eqn:dotLagCell}
\bx\cdot \z_i+p_i \geqslant \bx\cdot \z_j+p_j
\end{align}
for all $j\in\{1,\ldots,n\}$. The unknowns $p_j$ play the role of the weights used in Definition \ref{def:LagCells}.
(In the non-periodic case, the weights $w_i$ are related to the $p_i$ by $p_i = \frac 12  w_i - \frac 12 |\z_i|^2$, and
the modified geopotential $P$ is given by $P(\bx) = \max_i(\bx \cdot \z_i + p_i).$)
The task is to find $(p_1,\ldots,p_n)$ such that $C_i$ has given area 
$\overline{m}_i$
for each $i\in\{1,\ldots,n\}$. 
This is achieved by applying the nonlinear conjugate gradient method to the objective function
\[
\sum_{i=1}^n \left( |C_i| - \overline{m}_i \right)^2.
\]
At the first iteration, $p_i$ are defined such that all cells have positive area {(see \cite[Equation 5.10]{cullen2021mathematics})}. This initial guess is comparable to the initial guess for the weights used at the first time step in our algorithm. At each iteration, the cells $C_i$ are constructed using the divide-and-conquer algorithm of Preparata and Hong \cite{preparata1977convex}.

\section{Results}\label{sect:results}
{In this section we use the geometric method to simulate a shear instability and the formation of an atmospheric front.}
In what follows, the root mean square of the meridional velocity $v$ ($\mathrm{RMS}v$) at time $t$ is given by
\begin{align}
\label{eq:DefRMSv}
	{\mathrm{RMS}v}=
	\sqrt{\frac{1}{\vert\Omega\vert}\int_{\Omega}\left\vert v(\mathbf{x},t) \right\vert^2 \,\rd \mathbf{x}},
\end{align}
where $\vert \Omega \vert$ is the area of $\Omega$.

\subsection{Parameters and initial conditions}
\label{subsec:ParamValues}
In all simulations we used the following physical parameters, taken from \cite{visram2014framework}:
\begin{gather*}
	L=10^6\,\text{m}, 
	\qquad
	f=10^{-4}\,\text{s}^{-1},
	\qquad
	g=10\,\text{ms}^{-2}, 
	\qquad
	\theta_0=300\,\text{K},
	\\
	N=0.005\,\text{s}^{-1}  
	\qquad
	s=-3\times 10^{-6}\,\text{m}^{-1}\text{K},
	\qquad
	a=-7.5\, \text{ms}^{-1}.
\end{gather*}
We considered four initial conditions $\nabla P(\bx,0)$ of the form
\begin{align*}
    \nabla P(\bx,0) = \nabla \overline{P}(\bx)+G(\bx),
\end{align*}
{where $G$ is}
a small perturbation of the steady {shear flow}
\eqref{eq:Psteady}.
The perturbations $G$ are related to perturbations $\widetilde{\theta}$ of the potential temperature and $\widetilde{v}$ of the meridional velocity by
\begin{align}
    G(\bx)=
    \begin{pmatrix}
\frac{1}{f}\widetilde{v}(\mathbf{x})
\\
\frac{g}{f^2\theta_0} \widetilde{\theta}(\mathbf{x})
\end{pmatrix}.
\end{align}
In each case, a corresponding domain height $H$ was chosen, either in line with the linear instability analysis, or to enable comparison with the literature: see Table \ref{tab:ICs}.

The perturbations $G_{\mathrm{u}}$ and $G_{\mathrm{s}}$, listed in Table \ref{tab:ICs} and defined below, are normal modes of equation \eqref{eqn:SGESEul} linearised around the steady state \eqref{eq:Psteady}. Full details of the linearisation are contained in {Appendix \ref{Sec:LinearInstabilityAnalysis}.}
Define the \emph{Burger number} by
\[
\mathrm{Bu} = \frac{NH}{fL}.
\]
As shown in Section \ref{app:unstable}, exponentially growing normal modes exist only when the Burger number is less than a critical value $\mathrm{Bu}_{\mathrm{crit}}$, otherwise all normal modes are oscillatory. The domain height $H$ in Table \ref{tab:ICs} corresponding to $G_{\mathrm{u}}$ is chosen so that $\mathrm{Bu}< \mathrm{Bu}_{\mathrm{crit}}$, the maximum growth rate is achieved by the lowest frequency normal mode
$G_{\mathrm{u}}$, and all other normal modes are either decaying or oscillatory. Conversely, the domain height in 
Table \ref{tab:ICs} corresponding to $G_{\mathrm{s}}$ is chosen so that $\mathrm{Bu}>\mathrm{Bu}_{\mathrm{crit}}$,
and the lowest frequency normal mode $G_{\mathrm{s}}$ is oscillatory and has a wave speed of one domain length (i.e. $2L$) every $16$ days. 

\begin{table}
	\begin{center}
	\begin{tabular}{|c|c|c|c|}
		\hline
	 	Source                                                        & $G$                          & $H$ & $\mathrm{Bu}$ \\
		\hline
		Williams (1967) \cite{williams1967} 						  & $G_{\mathrm{u}}$ (see \eqref{eq:nablaP_u}) & $10224.85\,\mathrm{m}$  & $0.5112$ \\
		Appendix \ref{app:stable}                                     & $G_{\mathrm{s}}$ (see \eqref{eq:nablaP_s}) & $16374.56\,\mathrm{m}$  & $0.8187$ \\
		Visram \emph{et. al.} (2014) \cite{visram2014framework} 	  & $G_{\mathrm{V}}$ (see \eqref{eq:nablaP_V}) & $10^4\,\mathrm{m}$     & $0.5$ \\
		Cullen (2007) \cite[Equation 4.18]{cullen2007modelling} & $G_{\mathrm{C}}$ (see \eqref{eq:nablaP_C}) & $10^4\,\mathrm{m}$      & $0.5$ \\
		\hline
	\end{tabular}
	\caption{Initial conditions and corresponding parameter values used in this paper and in previous works.}
	\label{tab:ICs}
	\end{center}
\end{table}

Define the constants
\begin{align}
	\label{eqn:constA1}
	A_1 &:= \kappa \coth \kappa-1,\\
	\label{eqn:constA2}
	A_2 &:= \sigma(\kappa)
\end{align}
where
\begin{align}\label{eqn:kappa}
	\kappa & := \frac{\pi\mathrm{Bu}}{2},
\end{align}
and $\sigma:\R\to\R$ is given by
\begin{align}
\label{eqn:sigma}
	\sigma\left(\tilde{\kappa}\right):=
	\sqrt{
		\left| \left(\tilde{\kappa}-\tanh \tilde{\kappa} \right)  
		\left(\coth \tilde{\kappa}-\tilde{\kappa}\right)\right|
		}.
\end{align}
The unstable normal mode $G_{\mathrm{u}}$ is
\begin{align}\label{eq:nablaP_u}
    G_{\mathrm{u}}(\bx)=
    \begin{pmatrix}
\frac{1}{f}v_{\mathrm{u}}(\bx)
\\
\frac{g}{f^2\theta_0} \theta_{\mathrm{u}}(\bx)
\end{pmatrix},
\end{align}
where
\begin{align}\label{eq:thu}
	\theta_{\mathrm{u}}(\bx) = 
	\frac{aN\theta_0}{g}
	\left[
	A_1 \sinh \tfrac{\pi \mathrm{Bu}x_2}{H} \cos \tfrac{\pi x_1}{L} - A_2 \cosh \tfrac{\pi \mathrm{Bu}x_2}{H} \sin \tfrac{\pi x_1}{L} 	
	\right],
\end{align}
and
\begin{align*}
    v_{\mathrm{u}}(\bx)=
    -a \left[A_2 \sinh \tfrac{\pi\mathrm{Bu}x_2}{H}  \cos \tfrac{k \pi x_1}{L} + A_1 \cosh \tfrac{\pi\mathrm{Bu}x_2}{H} \sin \tfrac{k \pi x_1}{L} \right].
\end{align*}
The stable normal mode $G_{\mathrm{s}}$ is
\begin{align}\label{eq:nablaP_s}
G_{\mathrm{s}}(\bx)=
    \begin{pmatrix}
\frac{1}{f}v_{\mathrm{s}}(\mathbf{x})
\\
\frac{g}{f^2\theta_0} \theta_{\mathrm{s}}(\mathbf{x})
\end{pmatrix},
\end{align}
where
\begin{align*}
	\theta_{\text{s}}(\mathbf{x}) =
	\frac{aN\theta_0}{g}
	\cos \tfrac{\pi x_1}{L}
	\left[A_1 \sinh \tfrac{\pi \mathrm{Bu}x_2}{H}+ A_2 \cosh \tfrac{\pi \mathrm{Bu}x_2}{H}
	\right],
\end{align*}
and
\begin{align*}
	v_{\text{s}}(\mathbf{x}) =
	-a
	\sin \tfrac{\pi x_1}{L}
	\left[A_1 \cosh \tfrac{\pi \mathrm{Bu}x_2}{H}
	+ A_2 \sinh \tfrac{\pi \mathrm{Bu}x_2}{H}
	\right].
\end{align*}
The perturbation used in \cite{visram2014framework} and \cite{yamazaki2017vertical} is
\begin{align}\label{eq:nablaP_V}
    G_{\mathrm{V}}(\bx) =  G_{\mathrm{u}}(x_1,x_2/\pi).
\end{align}
Finally,
\begin{align}\label{eq:nablaP_C}
    G_{\mathrm{C}}(\bx)=
    \frac{gB}{\theta_0 f^2} \sin\left(\pi \left(\tfrac{x_1}{L}+\tfrac{x_2}{H}+\tfrac{1}{2}\right)\right)
    \begin{pmatrix}
    \tfrac{H}{L}
    \\
    1
    \end{pmatrix},
\end{align}
where
\begin{align*}
	B = 0.25 \, \text{K}.
\end{align*}
These, and the corresponding parameter values in Table \ref{tab:ICs}, are taken directly from the stated sources. Note that neither $G_{\mathrm{V}}$ nor $G_{\mathrm{C}}$ are normal modes of \eqref{eqn:SGESEul} linearised around the steady state \eqref{eq:Psteady}. They do, however, have the same horizontal wavelength as $\theta_u$.

Finally, in each simulation we specify three discretisation parameters: $n$, the number of cells in the spatial discretisation of $\Omega$; $h_{\mathrm{def}}$, the default time step size (see Algorithm \ref{alg:adaptiveTimeStep}); $\eta$, the percentage mass tolerance (see Algorithm \ref{alg:DN}).

\subsection{Unstable normal mode}
\label{Subsec:UNM}
The results reported in this subsection were obtained using the initial data given in Table \ref{tab:ICs}, Row 1, and the simulation parameters $\eta=0.01$, $h_{\mathrm{def}}=30$ and $n=2678$, unless otherwise stated. The dimensional growth rate of the unstable mode \eqref{eq:nablaP_u}, and the corresponding $\mathrm{RMS}v$, under the dynamics of the linearised equations \eqref{eq:AA7}-\eqref{eq:AA11} is
\begin{align*}
	\omega = -\frac{g s}{N \theta_0} \sigma(\kappa) = 0.53536 
	\text{ day}^{-1}.
\end{align*}
This is derived in the appendix: see equation \eqref{eq:FastestGrowthRate}.

\begin{figure}[t]
	\centering
    	\begin{subfigure}[t]{0.47\textwidth}
        	\includegraphics{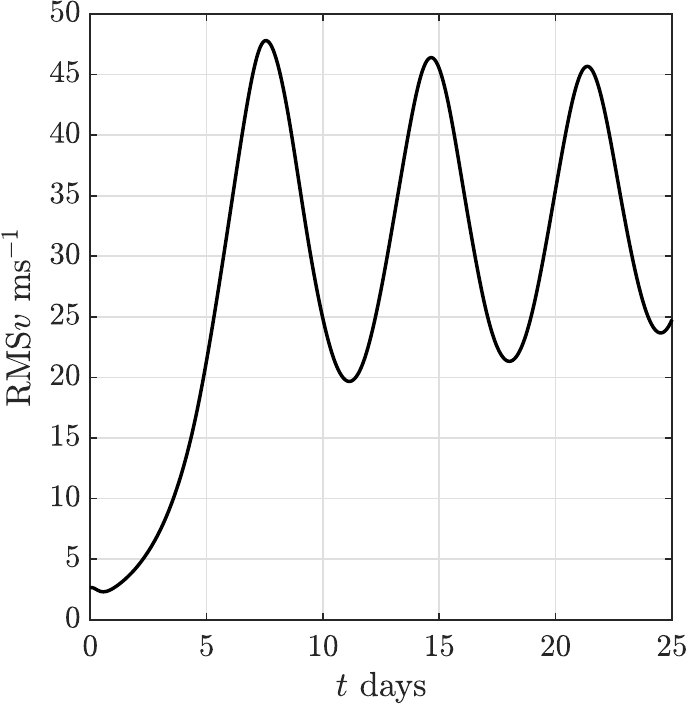}
        	\caption{Plot of the $\mathrm{RMS}v$ curve up to $t=25\,\mathrm{days}$. This figure illustrates the instability of the steady shear flow \eqref{eq:Psteady} in this parameter regime.
        	\label{fig:RMSvUnstable}}
        \end{subfigure}
        \hspace{0.02\textwidth}
        \begin{subfigure}[t]{0.47\textwidth}
            \includegraphics{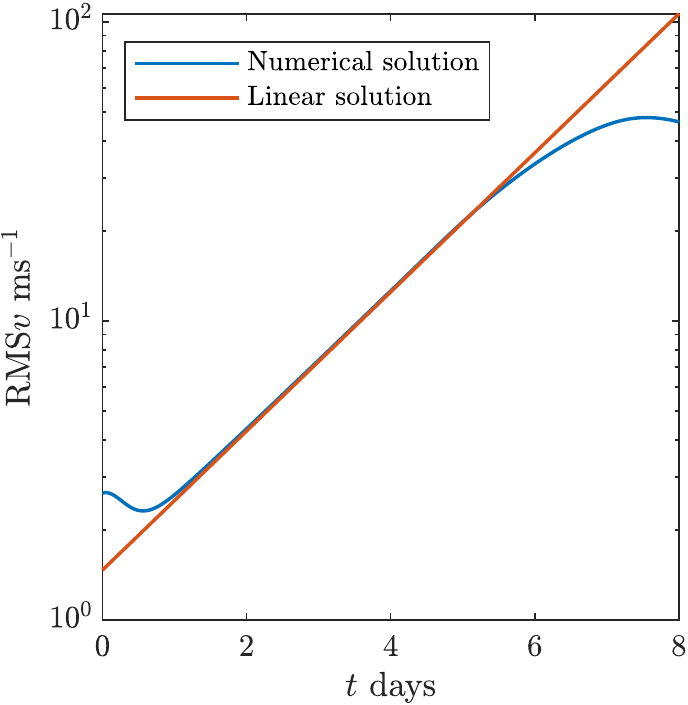}
        	\caption{Semi-log plot of the $\mathrm{RMS}v$ curve for the numerical solution of the nonlinear equation \eqref{eqn:SGESLagCts} (blue) and the exact solution of the linearised equations
        	{\eqref{eq:AA7}-\eqref{eq:AA11}}
        	(orange).
        	We see that the growth rate observed in simulations agrees well with the theoretical growth rate of the unstable perturbation.
        	}
        	\label{fig:RMSvGrowthRate}
        \end{subfigure}
    \caption{Plots of the $\mathrm{RMS}v$ curve for the numerical solution of \eqref{eqn:SGESLagCts} with initial data from Table \ref{tab:ICs}, Row 1, and simulation parameters $\eta=0.01$, $h_{\mathrm{def}}=30$ and $n=2678$. \label{fig:RMSv}}
\end{figure}

The time-evolution of the $\mathrm{RMS}v$ calculated from the simulation data is presented in Figure
\ref{fig:RMSvUnstable}.
Using equations \eqref{eqn:vThP} and \eqref{eq:nablaP}, the $\mathrm{RMS}v$ at time $t$ corresponding to a solution $\z=(\z_1,\ldots,\z_n)$ of \eqref{eqn:ODE} is given by
\begin{align}
    \mathrm{RMS}v=
    \sqrt{\frac{f^2}{|\Omega|}
    \sum_{i=1}^n
    {\int_{C_{i,\per}(t)}}
    \left\vert \Big(\bx-\z_i(t)-\bk_*\big(\bx,\z_i(t)\big)\Big)\cdot \mathbf{e}_1\right\vert^2\, \rd \bx},
\end{align}
where $\bk_*(\bx,\z_i(t))$ is defined by \eqref{eqn:kTranslate}, and $C_{i,\per}(t)$ is defined by \eqref{eq:brevity}.

After a short decline, the $\mathrm{RMS}v$ grows at the predicted rate $\omega$: see Figure \ref{fig:RMSvGrowthRate}. From around time $t=5$ days, this growth slows until a peak $\mathrm{RMS}v$ value is reached at time $t=7.5$ days. After a period of decline, the $\mathrm{RMS}v$ reaches a trough and proceeds to oscillate between the peak and trough values with a frequency of around 7 days. The initial decline in the $\mathrm{RMS}v$ is due to numerical errors incurred by the spatial discretisation: see Section \ref{sect:sensitivity}.

\begin{figure}
	\begin{subfigure}[b]{\textwidth}
		\centering
		\includegraphics[trim={0cm 0.5cm 1.75cm 0cm},clip,scale=1]{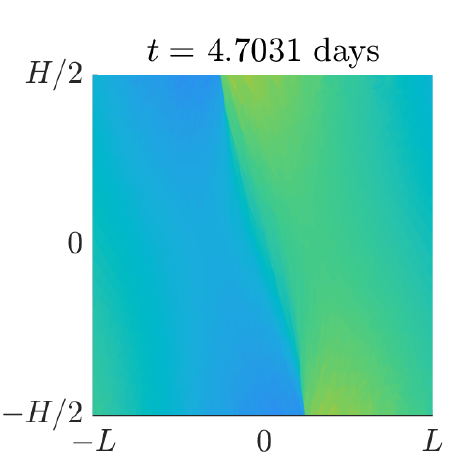}
		\includegraphics[trim={0cm 0.5cm 1.75cm 0cm},clip,scale=1]{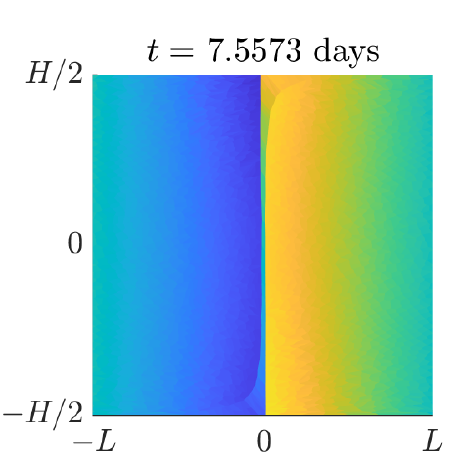}
		\includegraphics[trim={0cm 0.5cm 0cm 0cm},clip,scale=0.97]{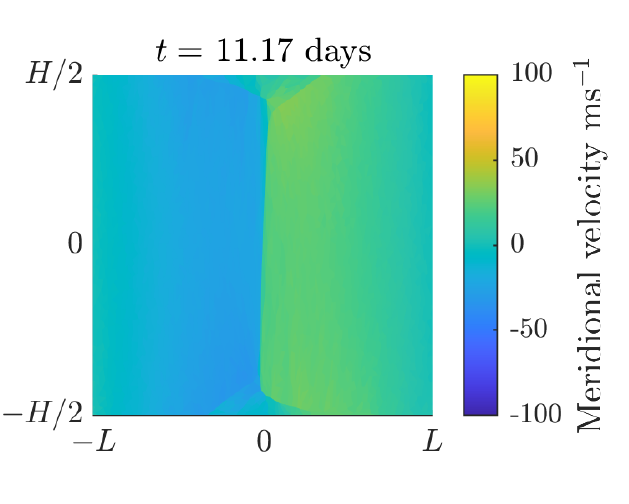}
		\vspace{0.3cm}
		\hspace{0.4cm}
		\includegraphics[trim={0cm 0cm 1.75cm 0cm},clip,scale=1]{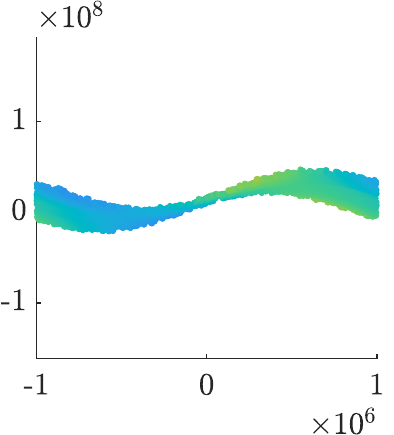}
		\hspace{0.5cm}
		\includegraphics[trim={0cm 0cm 1.75cm 0cm},clip,scale=1]{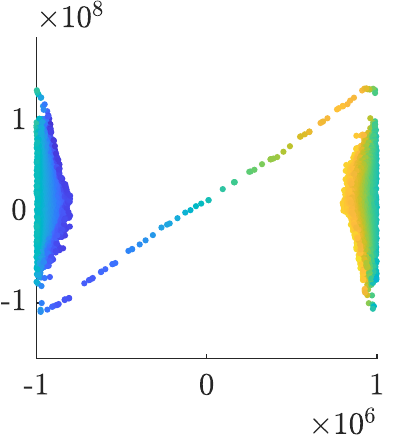}
		\hspace{0.5cm}
		\includegraphics[trim={0cm 0cm 0cm 0cm},clip,scale=1]{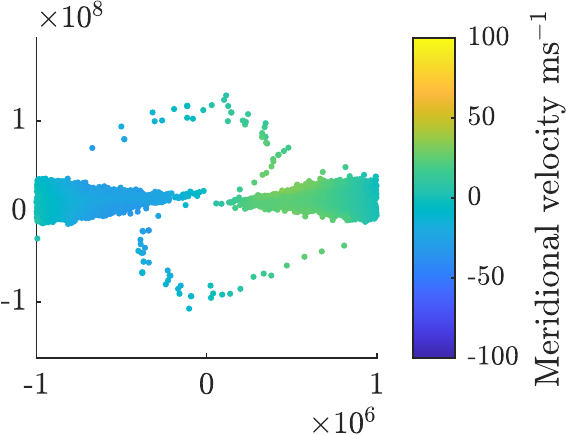}
	\end{subfigure}
	\begin{subfigure}[b]{\textwidth}
		\centering
		\includegraphics[trim={0cm 0.5cm 1.75cm 0cm},clip,scale=1]{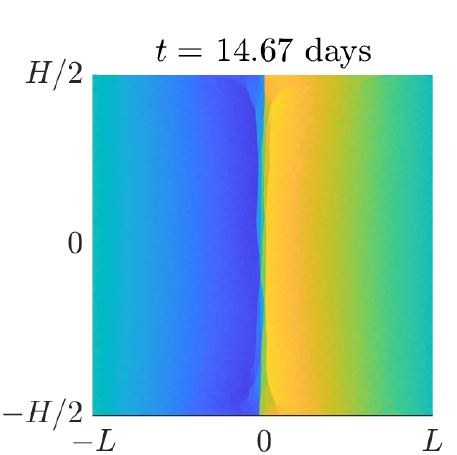}
		\includegraphics[trim={0cm 0.5cm 1.75cm 0cm},clip,scale=1]{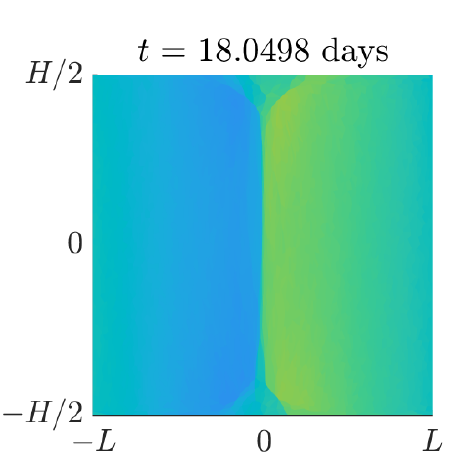}
		\includegraphics[trim={0cm 0.5cm 0cm 0cm},clip,scale=1]{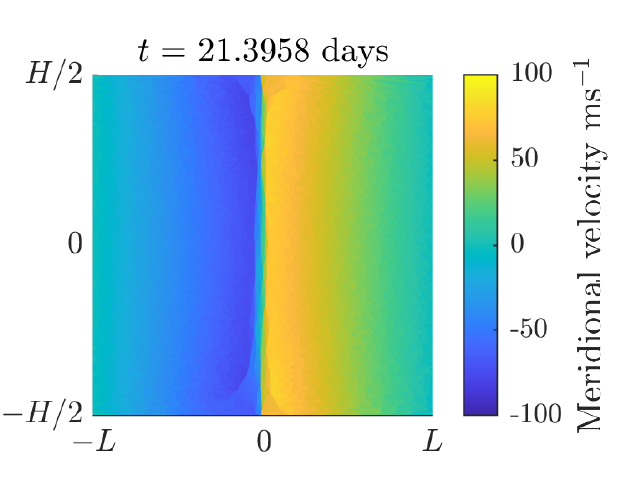}
		\vspace{0.3cm}
		\hspace{0.4cm}
		\includegraphics[trim={0cm 0cm 1.75cm 0cm},clip,scale=1]{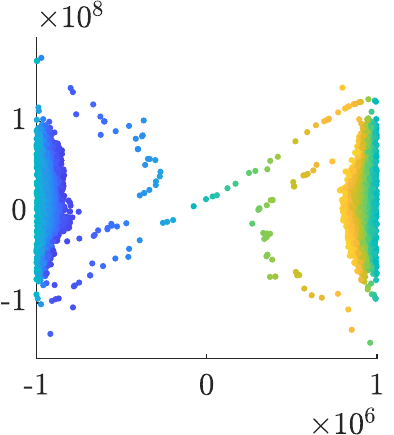}
		\hspace{0.5cm}
		\includegraphics[trim={0cm 0cm 1.75cm 0cm},clip,scale=1]{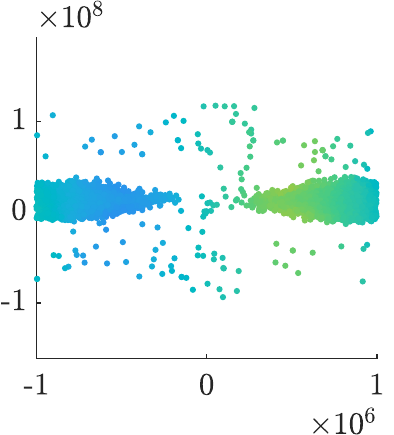}
		\hspace{0.5cm}
		\includegraphics[trim={0cm 0cm 0cm 0cm},clip,scale=1]{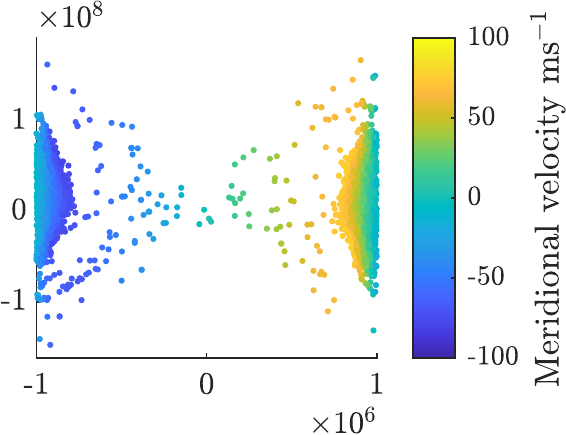}
	\end{subfigure}
	{\caption{Numerical solution of \eqref{eqn:SGESLagCts} with initial data from Table \ref{tab:ICs}, Row 1, and  simulation parameters $\eta=0.001$, $h_{\mathrm{def}}=30$ and $n=2678$. The Laguerre tessellations of the fluid domain are shown in rows 1 and 3, and the corresponding seeds $\z$ in geostrophic space are shown in rows 2 and 4. Both are coloured by the meridional velocity at the cell centroids. The first time $t=4.7031$ is halfway between $t=0$ and the time of the first peak of $\mathrm{RMS}v$. Subsequent times are chosen to coincide with peaks and troughs of the $\mathrm{RMS}v$. Observe the frontal discontinuity at time $t=7.5573$.}
	\label{fig:snapshotsUnstable}}
\end{figure}

The initial peak in the $\mathrm{RMS}v$ {curve
at time $t=7.5573$ days} coincides with the {formation}
of a frontal discontinuity, as seen in Figure \ref{fig:snapshotsUnstable}. 
Subsequent peaks/troughs of the $\mathrm{RMS}v$ occur when the frontal discontinuity is strongest/weakest, respectively.
Note that the meridional velocity $v$ displayed in Figure \ref{fig:snapshotsUnstable} is computed from the numerical solution $\z$ using \eqref{eqn:vThP} and \eqref{eq:nablaP}. It is a piecewise affine function with respect to the Laguerre cells $C_{i,\per}$,
which are shown in Figure \ref{fig:snapshotsUnstable}, Rows 1 and 3. They are coloured according to the meridional velocity at their centroids.

Up to first onset of frontogenesis, the distribution of seeds in geostrophic space appears to roughly approximate a two-dimensional subset of $\R^2$: see Figure \ref{fig:snapshotsUnstable}, $t=4.7031$ days. At the onset of frontogenesis, however, it is clear that a small subset of the seeds is distributed along a one-dimensional curve: see Figure \ref{fig:snapshotsUnstable}, $t=7.5573$ days. In other words, the frontal discontinuity in physical space appears to correspond to 
{a singular part of}
the \emph{potential vorticity}\footnote{The semi-geostrophic Eady Slice equation, and related semi-geostrophic systems, have been studied extensively in their potential vorticity formulation in the mathematical analysis literature: see, for example, \cite{benamou1998weak,loeper2006fully,feldman2013lagrangian,ambrosio2014global,feldman2017semi}. The relation of this viewpoint to the discrete formulation of the dynamics used in the geometric method is the subject of \cite{bourne2022semi}.} 
{measure} $\alpha$, which is defined 
{at time $t$ to be the push-forward measure $\alpha_t = \nabla P(\cdot,t) \# \mathcal{L}^2 \mres \Omega$}. 
For the case where $\alpha$ is non-singular (has a two-dimensional support), it is defined by
\begin{align*}
	\alpha_t(\nabla P(\mathbf{x},t)) = \det\left(D^2P(\mathbf{x},t)^{-1}\right).
\end{align*}

\begin{figure}[t]
	\centering
	\begin{subfigure}[t]{0.47\textwidth}
	    \includegraphics{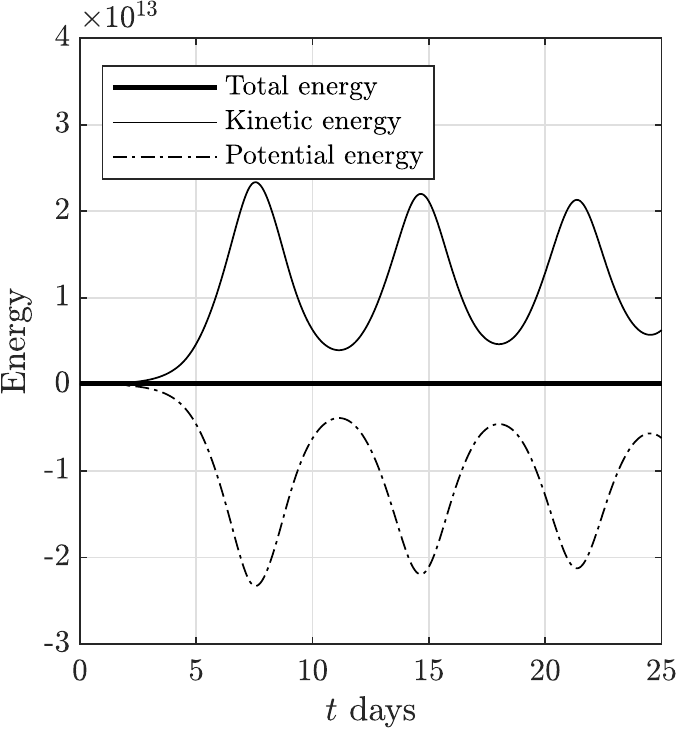}
	    \caption{Time-evolution of the total geostrophic energy \eqref{eqn:totalGeoEnergy}, kinetic energy \eqref{eqn:KE} and potential energy \eqref{eqn:PE} for the numerical solution of \eqref{eqn:SGESLagCts} with initial data from Table \ref{tab:ICs}, Row 1, and simulation parameters $\eta=0.001$, $h_{\mathrm{def}}=30$ and $n=2678$. This shows that our numerical scheme is energy-conserving to a high accuracy.}
	    \label{fig:energy1}
	\end{subfigure}
	\hspace{0.04\textwidth}
	\begin{subfigure}[t]{0.47\textwidth}
	    \includegraphics{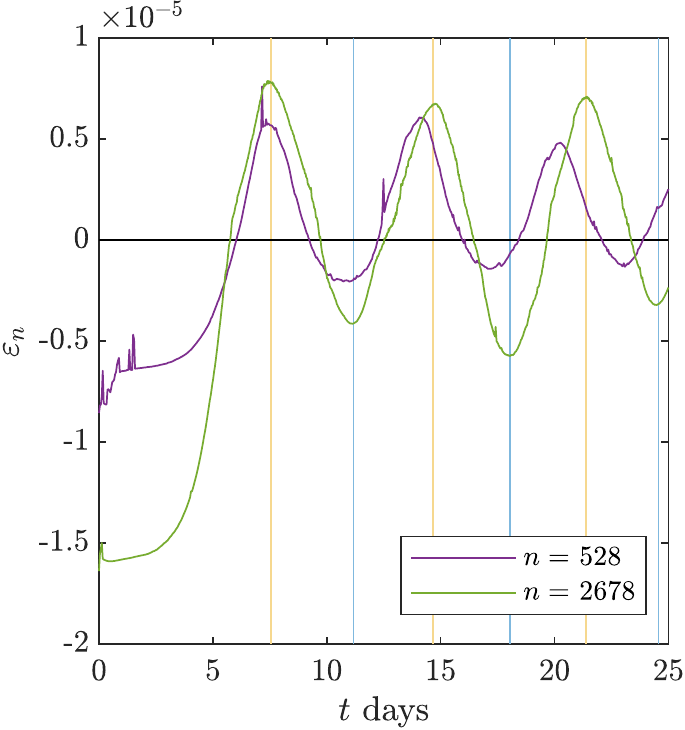}
	    \caption{Energy conservation error $\varepsilon_n$ (defined in \eqref{eq:epsilonn}) obtained using simulation parameters $\eta=0.001$ and $h_{\mathrm{def}}=30$. The $t$ coordinate of each orange or blue vertical line is a time when the $\mathrm{RMS}v$ curve corresponding to $n=2678$ is at a peak or trough, respectively. 
	We see that the energy conservation of the numerical method is worst at peaks and troughs of the $\mathrm{RMS}v$.}
	    \label{fig:energy2}
	\end{subfigure}
	\caption{Energy dynamics of our numerical solutions of \eqref{eqn:SGESLagCts}.}
	\label{fig:energy}
\end{figure}

The time evolution of the total geostrophic energy \eqref{eqn:totalGeoEnergy}, the total kinetic energy \eqref{eqn:KE} and the total potential energy \eqref{eqn:PE} is shown in Figure \ref{fig:energy1}. 
 Our results are coherent with the energy conservation property of the geometric method discussed in Section \ref{sect:Precovery}. To quantify to what extent the total geostrophic energy is conserved by 
 a numerical solution, we define the energy conservation error at time $t$ as
\begin{align}
\label{eq:epsilonn}
    \varepsilon_{n}(t)=\frac{\overline{\mathcal{E}}_n-\mathcal{E}_n(t)}{\overline{\mathcal{E}}_n},
\end{align}
where $\mathcal{E}_n(t)$ is the total geostrophic energy at time $t$ calculated from the numerical solution, $\overline{\mathcal{E}}_n$ is the temporal mean of this quantity, and $n$ is the number of seeds used in the simulation. 
For the values of $n$ listed in Table \ref{tab:sens},
$\vert\varepsilon_{n}(t)\vert<2\times 10^{-5}$ for all $t$.
Figure \ref{fig:energy2} demonstrates that $\vert\varepsilon_{n}(t)\vert$ is at a local maximum at times $t$ when the $\mathrm{RMS}v$ is at a peak or trough. By definition, peaks and troughs of the kinetic energy align with those of the $\mathrm{RMS}v$. Since the total energy is conserved, the total potential energy has the same behaviour as the total kinetic energy but with opposite sign, and energy is transferred from kinetic to potential, and vice-versa, over multiple lifecycles.


Numerical solutions of the Eady-Boussinesq vertical slice equations \eqref{eqn:EadySlice1} obtained in previous works using Eulerian \cite{yamazaki2017vertical} or semi-Lagrangian \cite{visram2014framework} methods exhibit significantly larger energy conservation errors. 
These losses in the total energy occur each time a frontal discontinuity forms. 
The energetic losses incurred in the results of \cite{yamazaki2017vertical} are attributed to the approximation of the advection of the meridional velocity $v$ in the numerical method. The geometric method has two advantages in this regard. First, it is a Lagrangian method so there is no need to approximate an advection term. 
Second, at each time step the Laguerre tessellation of the fluid domain is defined by the values of $\theta$ and $v$, and so it is adapted to the location of the front, regardless of the chosen resolution. This is reflected by the fact that the order of the energy conservation error was observed to be similar when using different numbers of seeds $n$.

The qualitative behaviour of our numerical solutions in \emph{physical} space coincide with those from previous works \cite{culrou,cullen2007modelling,visram2014framework,yamazaki2017vertical}. 
To \revA{the best of} our knowledge, the behaviour of the corresponding seeds in geostrophic space has not been previously \revA{illustrated}.

\subsection{Stable normal mode}
\label{Subsec:SNM}
The results reported in this subsection were obtained using the initial condition and physical parameters from row 2 of Table \ref{tab:ICs}, and the discretisation parameters $\eta =0.001$, $h_{\mathrm{def}}=30$ and $n=990$.


\begin{figure}
    \centering
		\includegraphics[trim={0cm 0cm 1.75cm 0cm},clip,scale=1]{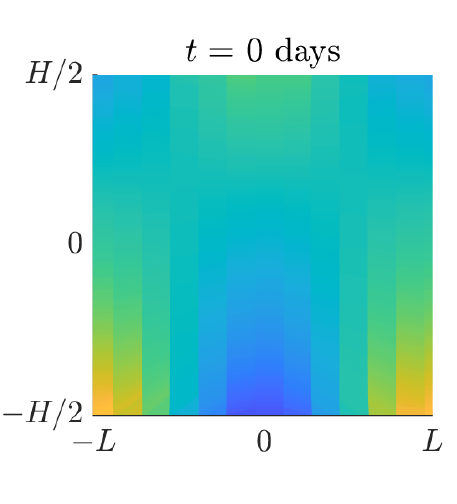}
		\includegraphics[trim={0cm 0cm 1.75cm 0cm},clip,scale=1]{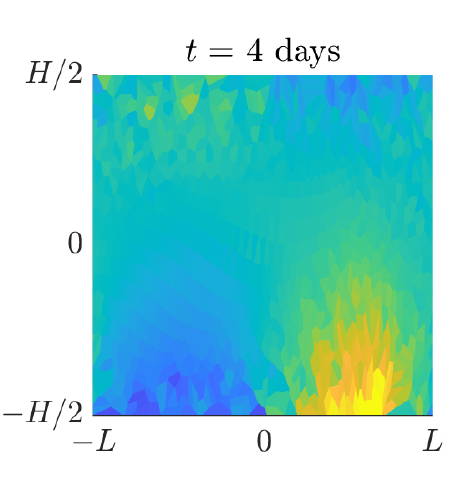}
		\includegraphics[scale=1]{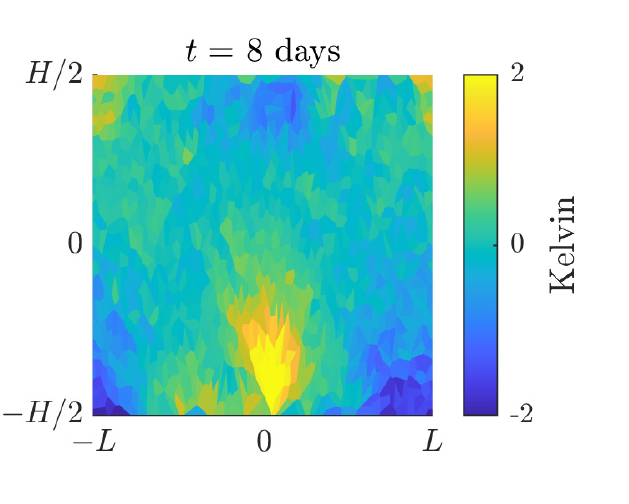}
	\caption{Plots of the potential temperature perturbation in the stable parameter regime with initial data from Table \ref{tab:ICs}, Row 2, and simulation parameters $\eta=0.001$, $h_{\mathrm{def}}=30$ and $n=990$. We see that the initial potential temperature perturbation propagates left with roughly constant amplitude. The numerically observed wave speed is in good agreement with the theoretical prediction of the linear instability analysis (Section \ref{app:stable}), which predicts that the wave takes $16$ days to cross the domain. In particular, the yellow patch in the bottom right corner at $t=0$ crosses half the domain by day $8$.}
	\label{fig:phaseLarge}
\end{figure}

Under the dynamics of the linearised equations \eqref{eq:AA7}-\eqref{eq:AA11}, disturbances of the steady shear flow \eqref{eq:Psteady} corresponding to the stable 
perturbation
\eqref{eq:nablaP_s} propagate west with constant wave speed
\begin{align*}
	c_{1} = \frac{|s|gL}{N\theta_0\pi}\sigma(\kappa),
\end{align*}
where $\kappa$ and $\sigma$ are defined by \eqref{eqn:kappa} and \eqref{eqn:sigma}: see Appendix \ref{app:stable}.
With our choice of physical parameters this equates to one domain length (i.e. $2L$) every $16$ days. This supports our numerical results, as can be seen in Figure \ref{fig:phaseLarge}, where we plot the perturbation of the steady potential temperature.
Indeed, Figure \ref{fig:phaseLarge} shows that the large scale 
disturbance initially located at the right-most boundary of $\Omega$ (see the yellow patch in the bottom right corner of Figure \ref{fig:phaseLarge}, $t=0$) moves left and
is centred around the line $x_1=L/2$ at day $4$
and the line $x_1=0$ at day $8$.
Of course, discretisation errors are incurred. One way in which these manifest is as small scale disturbances. 
As shown in equation \eqref{eq:stablewavespeed},
the maximum wave speed of such (large-wavenumber) disturbances under the dynamics of the linearised equations \eqref{eq:AA7}-\eqref{eq:AA11} is 
\begin{align*}
	c_{\infty} = \frac{\vert s\vert gH}{2f\theta_0},
\end{align*}
or approximately $0.3537$ domain lengths (i.e. $2L$) per day. While we do not include a corresponding figure, this can also be observed in our numerical results.

\subsection{Numerical evidence for the convergence of the Eady-Boussinesq vertical slice equations to the semi-geostrophic Eady slice equations}
\label{subsec:ComparisonLiterature}
We now compare our numerical solutions of the SG Eady slice equations \eqref{eqn:SGESLagCts} to the numerical solutions of the Eady-Boussinesq vertical slice equations \eqref{eqn:EadySlice1} obtained in \cite{nakamura1989nonlinear} and \cite{visram2014framework}.

As outlined in \cite[Section 2.3]{visram2014framework}, the SG  Eady slice equations \eqref{eqn:SGES1}-\eqref{eqn:vThRel0} can be understood formally as a small Rossby number approximation
of \eqref{eqn:EadySlice1}. 
This is made concrete by a rescaling argument using a scaling parameter $\beta$ such that the limit $\text{Ro}\to 0$ corresponds to the limit $\beta\to 0$. In \cite{visram2014framework}, numerical solutions of \eqref{eqn:EadySlice1} are obtained for a sequence of decreasing values of $\beta$ using a semi-Lagrangian method. 
The resulting $\mathrm{RMS}v$ curves are then compared to that from \cite[Figure 4.6]{cullen2007modelling}, which was obtained by solving the SG Eady slice equations \eqref{eqn:SGESLagCts} numerically using the geometric method. A similar programme
is followed in \cite{yamazaki2017vertical} using a compatible finite element method to solve \eqref{eqn:EadySlice1} 
numerically.
It is observed that the maximum amplitude of the $\mathrm{RMS}v$ obtained in \cite{cullen2007modelling} is much greater than that obtained in both subsequent studies, even for small $\beta$ and high-resolution simulations. This would not support the hypothesis that the SG Eady slice equations are
the limit of the Eady-Boussinesq vertical slice equations as $\text{Ro}\to 0$.
However, the reason for the discrepancy between the $\mathrm{RMS}v$ curves is that the physical parameters used in \cite{cullen2007modelling} are not the same as those used in the two subsequent papers, despite being reported as such.

\begin{figure}[t]
    \centering
    		\includegraphics[width=\textwidth,trim={0cm 1.5cm 0cm 3.65cm},clip]{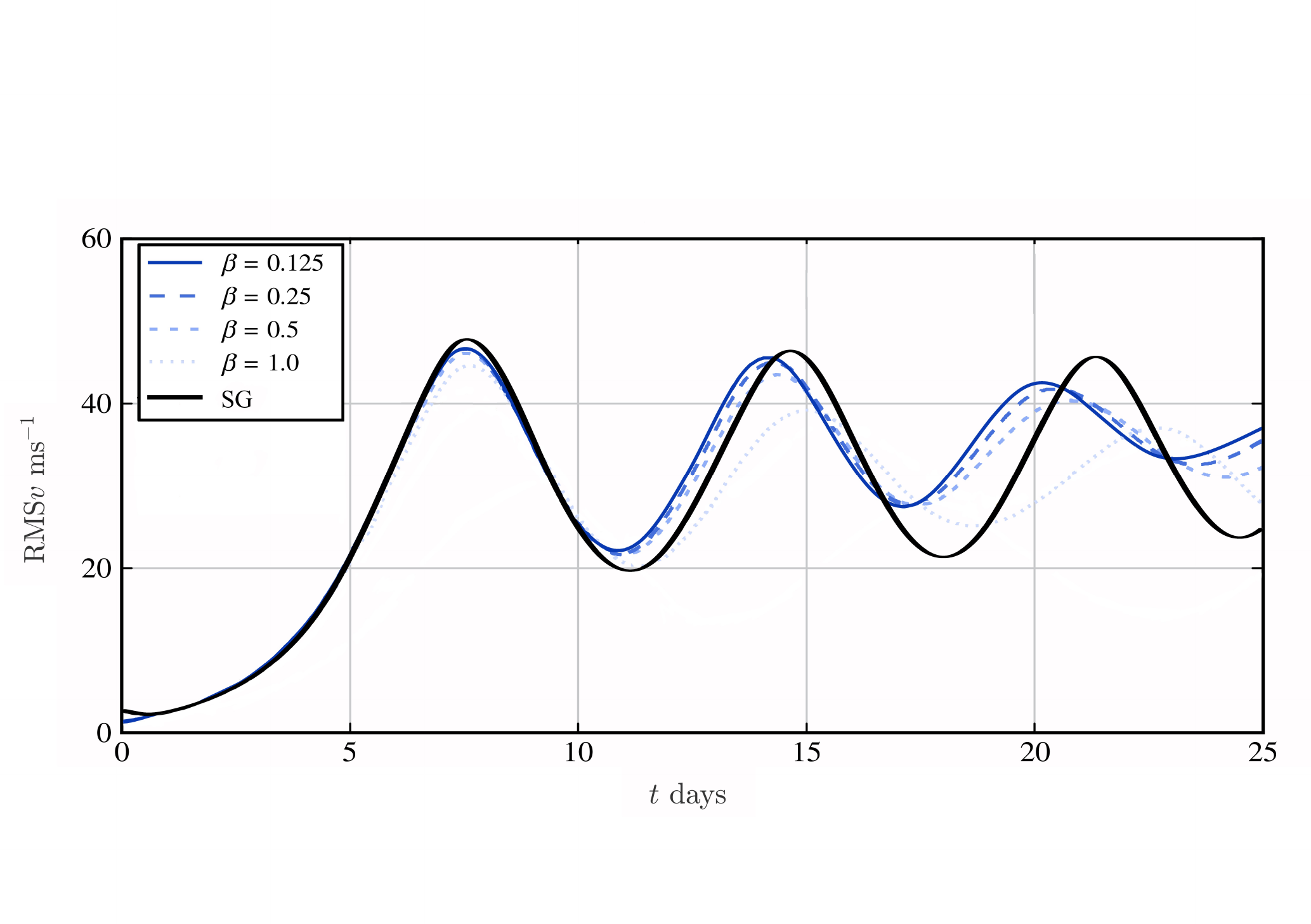}
    	{\caption{Comparison of the $\mathrm{RMS}v$ curve for the numerical solution of the SG Eady slice equations \eqref{eqn:SGESLagCts} (solid black curve) and the Eady-Boussineq vertical slice equations \eqref{eqn:EadySlice1} (blue curves, taken from \cite[Figure 5]{visram2014framework}). Equation \eqref{eqn:SGESLagCts} was solved with initial data from Table \ref{tab:ICs}, Row 1, and simulation parameters $\eta=0.01$, $h_{\mathrm{def}}=30$ and $n=2678$. 
    	The limit $\beta\to 0$ corresponds to the large scale limit $\mathrm{Ro}\to 0$. As $\beta \to 0$, the blue curves tend towards the solid black curve, at least for small times. This supports the hypothesis that the SG Eady slice equations are the small Rossby number limit of the Eady-Boussinesq vertical slice equations.
    	}
    	\label{fig:RMSvComparisonNormalMode}}
\end{figure}

To correct the comparisons made in \cite{visram2014framework} and \cite{yamazaki2017vertical}, we implement the geometric method using the physical parameters listed in those papers. Note that the initial condition used in both \cite{visram2014framework} and \cite{yamazaki2017vertical} is that of Table \ref{tab:ICs}, Row 3, not the most unstable mode discussed in Section \ref{Subsec:UNM}. When using this initial condition, it is therefore some time before the fastest growing unstable normal mode dominates the numerical solution and the expected growth rate of the $\mathrm{RMS}v$ is achieved. To account for this, in \cite{visram2014framework} and \cite{yamazaki2017vertical} the numerical solutions are translated backwards in time so that the maximum value of $v$ at $t=0$ days matches that of \cite{nakamura1989nonlinear}. We use the normal mode initial condition from Table \ref{tab:ICs}, Row 1, and do not translate in time.
The discretisation parameters we use are
$\eta = 0.01$, $h_{\mathrm{def}} = 30$ and $n= 2678$.
The resulting $\mathrm{RMS}v$ curve is plotted against those from \cite{visram2014framework} in Figure \ref{fig:RMSvComparisonNormalMode}\footnote{
Note that a similar figure appears in \cite{cullen2021mathematics}. There, however, the $\mathrm{RMS}v$ curve comes from the numerical solution of the SG Eady slice equations \eqref{eqn:SGESLagCts} obtained using our implementation with the initial condition from Table \ref{tab:ICs}, Row 4, and discretisation parameters $\eta = 0.01$, $h_{\mathrm{def}} = 15$ and $n= 1575$. Since the initial condition is not the most unstable normal mode, for comparison, the $\mathrm{RMS}v$ curve is translated backwards in time as in \cite{visram2014framework} and \cite{yamazaki2017vertical}.
}. We see that as $\beta \to 0$, the $\mathrm{RMS}v$ curves of \cite{visram2014framework} tend towards 
the $\mathrm{RMS}v$ curve of our numerical solution of the SG Eady slice equations \eqref{eqn:SGESLagCts} (solid black curve in Figure \ref{fig:RMSvComparisonNormalMode}), at least for small times. This supports the hypothesis that \eqref{eqn:SGES1}-\eqref{eqn:vThRel0} are the small Rossby number limit of \eqref{eqn:EadySlice1}. In particular, except at very early times where the SG solution is affected by discretisation errors, the initial growth rate of the $\mathrm{RMS}v$ curves from \cite{visram2014framework} approaches that of the SG numerical solution as $\beta\to 0$. This supports the initialisation procedure employed in \cite{visram2014framework} and \cite{yamazaki2017vertical}. For a fair and more detailed comparison, however, it would be desirable to perform further experiments using the semi-Lagrangian method of \cite{visram2014framework} and the compatible finite-element method of \cite{yamazaki2017vertical} with the normal mode initial condition so that every simulation uses the same initial condition and there is no need to time-translate the solutions.

\subsection{Sensitivity analysis}
\label{sect:sensitivity}
In this section we investigate the sensitivity of our numerical solution 
to the simulation parameters. These are $\eta$, the percentage area tolerance supplied to the semi-discrete optimal transport solver; $h_{\mathrm{def}}$, the default time step size in seconds; $n$, the number of seeds. The initial condition and physical parameters used to generate the results in this section are those from Table \ref{tab:ICs}, Row 1.

\begin{table}[t]
\begin{center}
\resizebox{\columnwidth}{!}{
	\begin{tabular}{|c c |c c c c|c c c c|}
		\hline
		 &             &   \multicolumn{4}{c|} {Square root of Sinkhorn loss}                                &                       
		 \multicolumn{4}{c|}{Weighted Euclidean distance} 
		 \\    
		 &             & $t=2\,\mathrm{d}$  & $t=4\,\mathrm{d}$ & $t=6\,\mathrm{d}$ & $t=8\,\mathrm{d}$ & $t=2\,\mathrm{d}$  & $t=4\,\mathrm{d}$ & $t=6\,\mathrm{d}$ & $t=8\,\mathrm{d}$ \\
		\hline
		
		$\eta$ & $1$     & $4.93 e{-4}$  & $7.41 e{-4}$ & $2.34 e{-3}$ & $5.49 e{-3}$
		
		                 & $9.54 e{-4}$  & $2.96 e{-3}$ & $2.52 e{-2}$ & $5.67 e{-2}$ \\
		                 
		   	   & $0.1$   & $4.16 e{-6}$ & $2.05 e{-4}$  & $1.44 e{-3}$ & $2.99 e{-3}$
		   	   
		   	             & $1.87 e{-5}$ & $1.10 e{-3}$ & $1.30 e{-2}$ & $2.96 e{-2}$ \\
		   	             
	           & $0.01$  & $9.83 e{-7}$ & $9.45 e{-5}$  & $1.15 e{-3}$ & $2.72 e{-3}$
	           
	                     & $7.85 e{-6}$ & $6.81 e{-4}$ & $1.08 e{-2}$ & $2.75 e{-2}$\\
	                     
	           & $0.001$  
	           & 		-	    & 		-      &  		-     & - 
	                     & 		-	    & 		-      & 		-     & - \\
          		\hline
		$h_{\mathrm{def}}\, \mathrm{(s)}$& $60$ & $8.59 e{-6}$  & $7.91 e{-3}$ & $1.52 e{-3}$ & $2.27 e{-3}$   
		
		            				 & $1.91 e{-5}$ & $1.69 e{-2}$ & $1.04 e{-2}$ & $2.68 e{-2}$ \\ 

	         			  & $30$ & $1.52 e{-6}$ & $7.48 e{-5}$ & $2.74 e{-3}$ & $2.49 e{-3}$   

	                		      	 & $6.64 e{-6}$ & $1.00 e{-3}$ & $1.28 e{-2}$ & $3.00 e{-2}$ \\ 

			 			  & $15$ & $1.64 e{-6}$ & $4.46 e{-5}$ & $1.08 e{-3}$ & $3.15 e{-3}$

			                      & $6.50 e{-6}$ & $6.22 e{-4}$ & $8.15 e{-3}$ & $2.43 e{-2}$\\

			              & $7.5$ &  - & - &- &- &- &- &- &- \\
					\hline
		$n$ & $528$  & $8.74 e{-3}$& $6.55 e{-3}$ & $6.95 e{-3}$ & $1.08 e{-2}$    
		             &    -        &      -       &      -       &      -       \\ 
			& $944$  & $4.12 e{-3}$& $3.71 e{-3}$ & $4.78 e{-3}$ & $8.00 e{-3}$    
			         &    -        &      -       &      -       &      -       \\ 
			& $1470$ & $2.20 e{-3}$& $2.59 e{-3}$ & $2.97 e{-3}$ & $3.33 e{-3}$
			         &    -        &      -       &      -       &      -       \\
			& $2124$ & $1.41 e{-3}$& $2.07 e{-3}$ & $3.29 e{-3}$ & $3.97 e{-3}$ 
			         &    -        &      -       &      -       &      -       \\
	        & $2678$ & - & - & - & - & - & - & - & - \\
		\hline
	\end{tabular}
	}
	\captionof{table}{
	Relative errors of the numerical solution of \eqref{eqn:SGESLagCts}
	with initial data from Table \ref{tab:ICs}, Row 1 for different values of the discretisation parameters $\eta$ (optimal transport tolerance), $h_{\mathrm{def}}$ (default time step) and $n$ (number of seeds). The second column gives the relative errors at different time steps with respect to the Sinkhorn loss (defined in \cite[Theorem 1]{genevay2018learning}) with regularisation parameter $\varepsilon=0.001$. The third column gives the relative error
	with respect to the weighted Euclidean norm (defined in  \eqref{eq:WEN}). 
	The errors are normalised by the factor in \eqref{eq:NormFactor}.
	In Row 2, $h_{\mathrm{def}}=30$ and $n=1470$; in Row 3 $\eta=0.001$ and $n=1470$; in Row 4 $\eta=0.01$ and $h_{\mathrm{def}}=30$.
	Numerical solutions corresponding to each row were compared to the numerical solution obtained using the finest discretisation parameter in that row.
	There are no values in the bottom right section of the table since the weighted Euclidean distance can only be computed between two numerical solutions with the same value of $n$. The relative errors are small and, in most cases, decrease as $\eta$ and $h_{\mathrm{def}}$ decrease and as $n$  increases. They are plotted in Fig.~\ref{fig:relativeErrors}.
	}
    \label{tab:sens}
\end{center}
\end{table}

\begin{figure}[t]
		\centering
		\includegraphics[width=\textwidth]{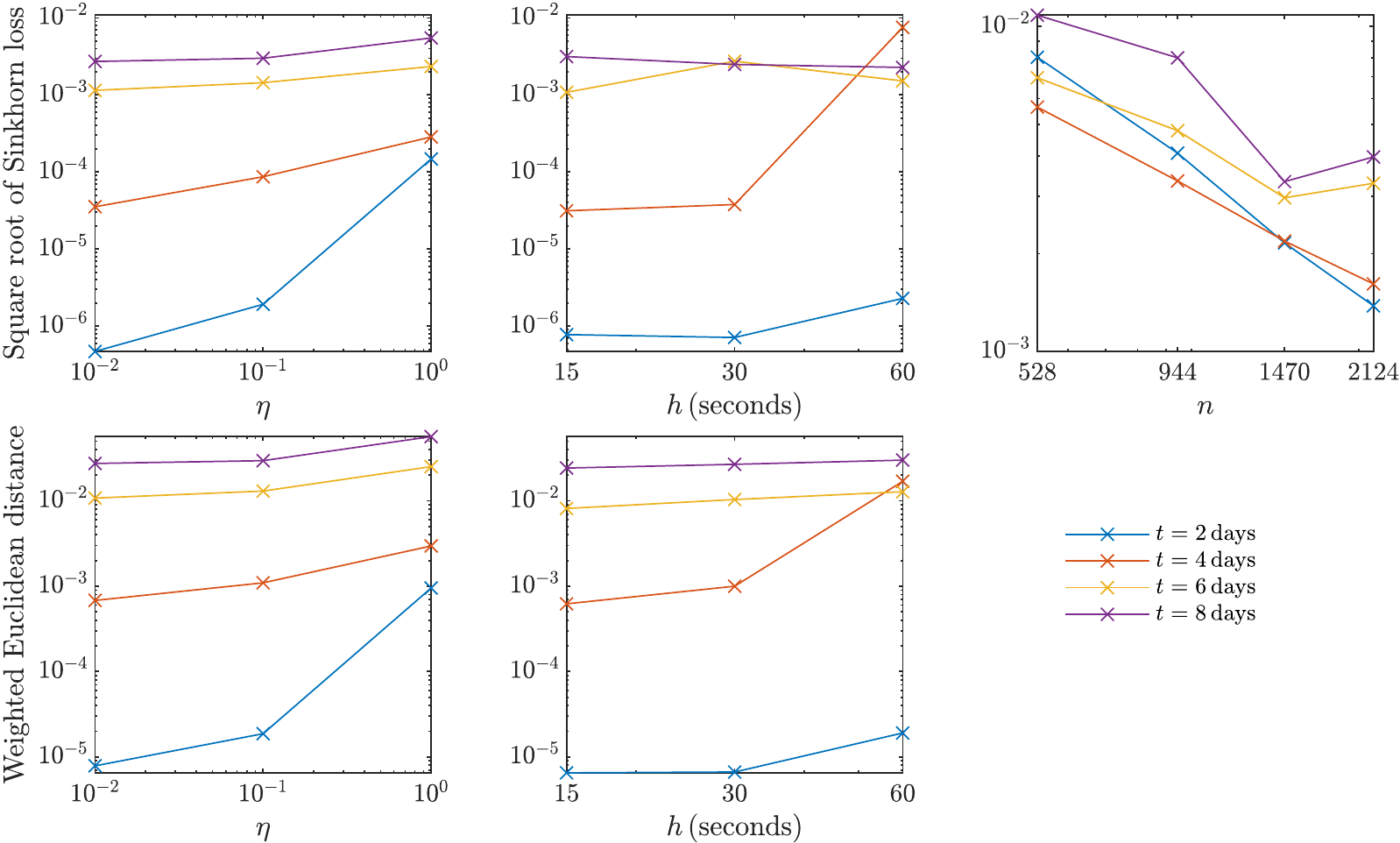}
	{\caption{Plots of the relative errors from Table \ref{tab:sens}.}
	\label{fig:relativeErrors}
	}
\end{figure}

\paragraph{Accuracy tests.} To evaluate the accuracy of our numerical solutions of \eqref{eqn:SGESLagCts}, we ran three sets of simulations, varying each parameter in turn. For the first set of simulations in which $\eta$ was varied, we used $h_{\mathrm{def}}=30$, and $n=1470$; in the second set of simulations in which $h_{\mathrm{def}}$ was varied, we used $\eta=0.001$, and $n=1470$; in the third set of simulations in which $n$ was varied, we used $\eta=0.01$ and $h_{\mathrm{def}}=30$. The remaining parameter values are listed in the first column of Table \ref{tab:sens}.
For each set of simulations, the numerical solutions were compared to that obtained using the discretisation parameter giving the finest discretisation, which is listed in the final row of the corresponding section of Table \ref{tab:sens}.

Next we state how we compare solutions with different discretisation parameters.
Associated to each solution $\z$ of the ODE \eqref{eqn:ODE} with corresponding target masses $\obm=(\om_1,\ldots,\om_n)$ is a time-dependent family of discrete distributions
\begin{align*}
\sum_{i=1}^n\om_i\delta_{\z_i(t)},
\end{align*}
known as the \emph{potential vorticity}. A natural way to quantify the discrepancy between two discrete distributions 
\begin{align*}
\sum_{i=1}^nm_i\delta_{\z_i}
\quad 
\text{and}
\quad
\sum_{j=1}^{n^{\prime}}m^{\prime}_j\delta_{\z_j^{\prime}}
\end{align*}
is by using the 
\revB{Wasserstein $2$-distance (see, for example, \cite[Chapter 5]{santambrogio2015optimal}). The \emph{Sinkhorn loss} introduced in \cite[Theorem 1]{genevay2018learning} provides an approximation of the Wasserstein $2$-distance that can be computed easily and quickly by solving a regularised optimal transport problem. Indeed, its square-root converges to the Wasserstein $2$-distance as the strength of the regularisation goes to zero \cite{feydy2019interpolating}.} If $n=n^{\prime}$ and $m_i=m_i^{\prime}$ for all $i\in\{1,\ldots,n\}$, one can also consider the weighted Euclidean distance between $\z$ and $\z^{\prime}$ given by
\begin{align}
\label{eq:WEN}
\sqrt{\sum_{i=1}^nm_i\left\vert\z_i-\z_i^{\prime}\right\vert^2},
\end{align}
which provides an upper bound on the Wasserstein 2-distance and is much simpler to compute. 
For each set of simulations, \revB{we computed both the square root of the Sinkhorn loss with regularisation parameter $\varepsilon=0.001$ and the weighted Euclidean distance} at several times $t$ and normalised them by the
factor
\begin{align}
\label{eq:NormFactor}
	\left(|\Omega|\underset{i}{\max}\|\z_i(t)\|^2_{\infty}\right)^{-\frac 12},
\end{align}
where $\z$ is the numerical solution of \eqref{eqn:ODE} obtained using the simulation parameter giving the finest discretisation, and $\|\z_i(t)\|_{\infty}$ is the maximum absolute value of the components of $\z_i(t)$. These values are presented in Table \ref{tab:sens} and Figure \ref{fig:relativeErrors} and, because of the choice of normalisation, can be interpreted as relative errors.

The relative error decreases with $\eta$, and tends to increase over time, as expected (see Figure \ref{fig:relativeErrors}, top left and bottom left). When quantified using the weighted Euclidean distance, the relative error also decreases with the default time step size $h_{\mathrm{def}}$, but when quantified using the Sinkhorn loss, this is not true at later times (see Figure \ref{fig:relativeErrors}, top middle and bottom middle).
The relative error decreases with $n$ at $t=2\,\mathrm{days}$ and $t=4\,\mathrm{days}$, but appears to increase slightly for large values of $n$ at $t=6\,\mathrm{days}$ and $t=8\,\mathrm{days}$
(see Figure \ref{fig:relativeErrors}, top right).
Nevertheless, in all cases, the relative errors are small. The unexpected trends may occur because the reference solutions are numerical rather than exact, or because the time step sizes chosen by the adaptive method depend on the simulation parameters.

\begin{figure}
	\centering
	\begin{subfigure}[t]{0.48\textwidth}
	    \includegraphics[scale=1]{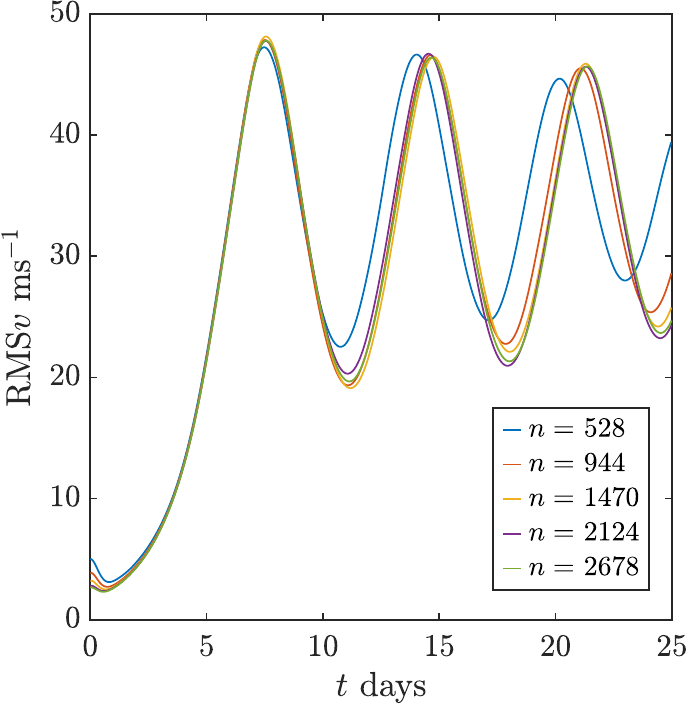}
	    \caption{$\mathrm{RMS}v$ curves for different numbers of seeds $n$ over multiple life cycles. The $\mathrm{RMS}v$ curves converge as $n$ increases.}
	\label{fig:RMSvSensitivity1}
	\end{subfigure}
	\hspace{0.02\textwidth}
	\begin{subfigure}[t]{0.48\textwidth}
	    \includegraphics[scale=1]{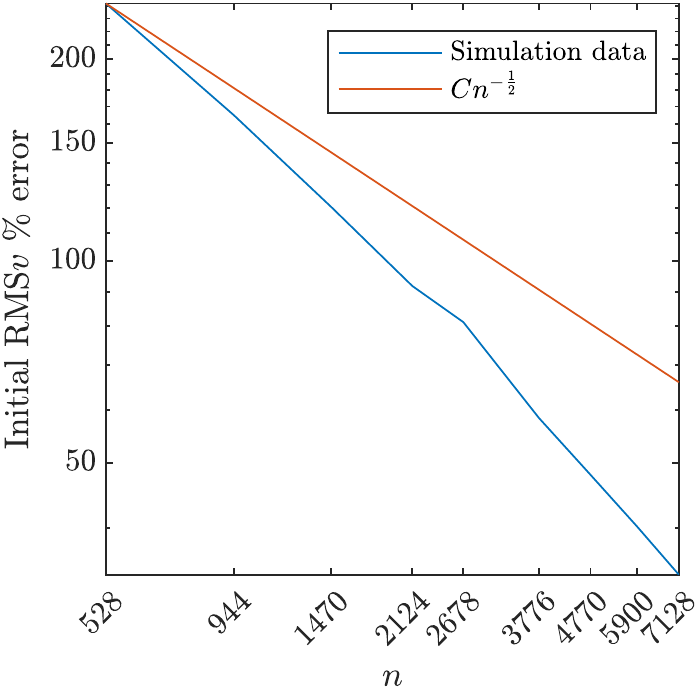}
	    \caption{Error of discretising the initial condition: see Section \ref{subsec:quantizeIC}. This log-log plot shows the $\mathrm{RMS}v$ percentage error at $t=0$ against the number of seeds $n$. {$C$ is a positive constant.} The discretisation error decays at a rate faster than $n^{-\frac 12}$.}
	\label{fig:RMSvSensitivity2}
	\end{subfigure}
	
	\begin{subfigure}[]{0.49\textwidth}
    \includegraphics{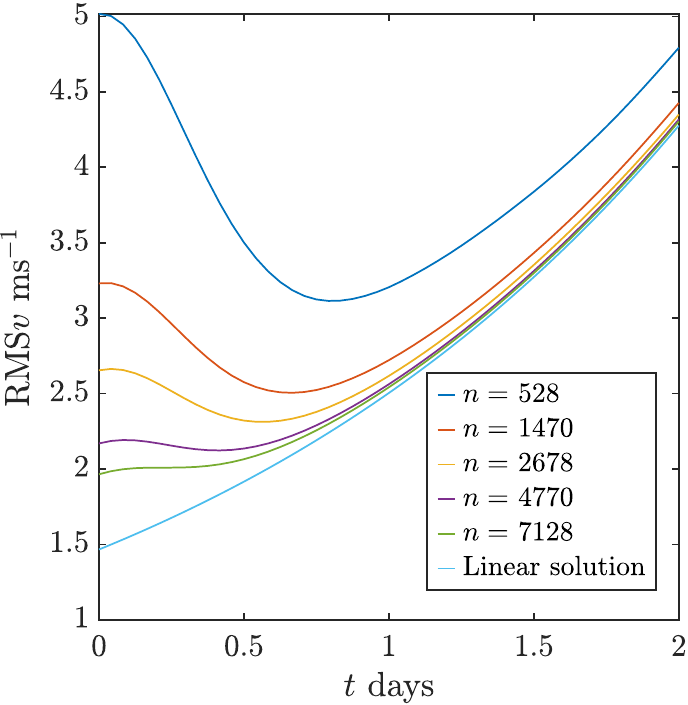}
    \caption{$\mathrm{RMS}v$ curves up to $t=2$ days with large $n$. The light-blue curve is the $\mathrm{RMS}v$ curve of the exact solution of the linearised equations \eqref{eq:AA7}-\eqref{eq:AA11}. We see that the initial dip in the $\mathrm{RMS}v$ is an artefact of the discretisation.}
    \label{fig:RMSvSensitivity3}
    \end{subfigure}
	\caption{Sensitivity of the $\mathrm{RMS}v$ to the number of seeds $n$. The initial data was taken from Table \ref{tab:ICs}, Row 1, and the discretisation parameters were $\eta=0.01$ and $h_{\mathrm{def}}=30$.}
	\label{fig:RMSvSensitivity}
\end{figure}

\paragraph{Behaviour of RMSv.}
Figure \ref{fig:RMSvSensitivity1} demonstrates that for sufficiently large $n$ the $\mathrm{RMS}v$ is well approximated over multiple life-cycles by our numerical results. We deduce from Figures \ref{fig:RMSvSensitivity2} and \ref{fig:RMSvSensitivity3} that the initial decline in the $\mathrm{RMS}v$ is due to errors incurred by the spatial discretisation of the initial condition (see Section \ref{subsec:quantizeIC}), which decrease at a rate faster than 
$n^{-\frac{1}{2}}$,
and do not significantly affect the behaviour of the $\mathrm{RMS}v$ after $t=2\,\mathrm{days}$.

\begin{figure}
	\centering
	\includegraphics[scale=1]{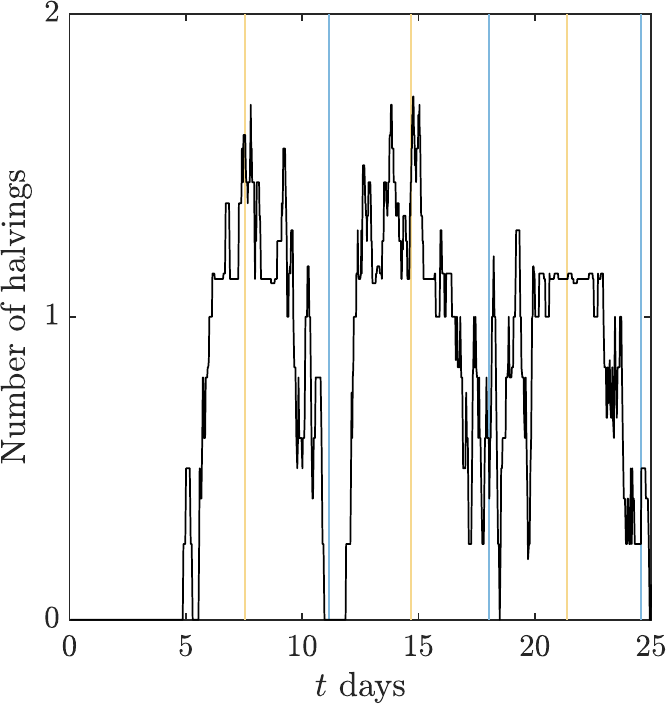}
	\caption{$4$-hour average number of time step halvings {per time step (see Step 2 in Algorithm \ref{alg:adaptiveTimeStep})} using a default time step size of $h_{\mathrm{def}}=30$ with $\eta=0.01$, $n=2678$
{and initial data from Table \ref{tab:ICs}, Row 1.}
	The $t$ coordinate of each orange or blue vertical line is a time when the corresponding $\mathrm{RMS}v$ curve is at a peak or trough, respectively. 
	{The number of time step refinements is maximum at peaks of the $\mathrm{RMS}v$ curve, which correspond to peaks in the strength of the frontal discontinuity. Therefore the adaptive time stepping scheme successfully identifies the frontal discontinuities.}
	}
	\label{fig:halvings}
\end{figure}

\paragraph{Adaptive time step size.} 
Figure \ref{fig:halvings} demonstrates that the number of times that the time step size is halved 
per iteration
by the adaptive time stepping algorithm correlates with the magnitude of the $\mathrm{RMS}v$, and hence the strength of the frontal discontinuity. This is to be expected because the vertical component of $\dot{\z}_i(t)$ is proportional to
{the meridional velocity}
$v(\bc_i(\z(t)),t)$ so, roughly speaking, as the $\mathrm{RMS}v$ increases, so does the magnitude of $\dot{\z}_i(t)$. 
As such, the first-order Taylor expansion of $\w_*$ becomes less accurate as the $\mathrm{RMS}v$ increases, so a smaller time step is needed in order for it to generate a good initial guess for the weights for Algorithm \ref{alg:DN}.

\section{Conclusion}
In this paper, we recast the geometric method for solving the SG Eady slice equations in the language of semi-discrete optimal transport theory (Section \ref{sect:geomMethod}), and develop a new implementation using the latest results from semi-discrete optimal transport theory and a novel adaptive time-stepping algorithm that is tailored to the ODE \eqref{eqn:ODE} (Algorithm \ref{alg:adaptiveTimeStep}).
The numerical solutions that we obtain via our implementation support the conjecture that weak solutions of the Eady-Boussinesq vertical slice equations \eqref{eqn:EadySlice1} converge to weak solutions of the semi-geostrophic Eady slice equations \eqref{eqn:SGES1}, \eqref{eqn:vThRel0} as $\mathrm{Ro} \to 0$ (Section \ref{subsec:ComparisonLiterature}).
Rigorous numerical tests in Section \ref{sect:sensitivity} demonstrate the sensitivity of the algorithm with respect to the discretisation parameters.
To clarify the use of different initial conditions in the literature on the Eady slice problem \cite{cullen2007modelling, williams1967, visram2014framework, yamazaki2017vertical}, we include a linear instability analysis of the steady shear flow \eqref{eq:SteadyState}, validating and extending the work of Eady \cite{eady1949long} (see Appendix \ref{Sec:LinearInstabilityAnalysis}). We perform simulations
in different physical parameter regimes, which verify the linear instability analysis (Sections \ref{Subsec:UNM} and \ref{Subsec:SNM}). 

The linear instability analysis provides benchmark initial conditions in both stable and unstable parameter regimes. Along with our implementation of the geometric method, this could be used in future work to carry out a more rigorous numerical study of the convergence of weak solutions of the Eady-Boussinesq vertical slice equations to weak solutions of the SG Eady slice equations as $\mathrm{Ro}\to 0$, and to explore the behaviour of solutions in different physical parameter regimes.

\appendix
\section{Derivation of the semi-discrete transport problem}
\label{app:OTderivation}
In this section we derive the semi-discrete optimal transport problem (Definition \ref{def:SDOT}) from the stability principle. 
Consider a geopotential $\phi:\R\times [-H/2,H/2]\to\R$ that is $2L$-periodic in the $x_1$ direction.
The modified geopotential $P:\R\times [-H/2,H/2]\to\R$ defined by
\begin{align*}
P(\bx)=\frac 12 x_1^2 + \frac 1f \phi(\bx)
\end{align*}
satisfies
\begin{align*}
\nabla P(\bx+\mathbf{k}) = \nabla P (\bx) + \mathbf{k}
\end{align*}
for all $\mathbf{k}\in K$ and all $\bx\in\R\times [-H/2,H/2]$ where $P$ is differentiable. If, in addition, $P$ is piecewise affine, then there exist $\z_i\in [-L,L)\times\R,\,i\in\{1,\ldots,n\}$, and a collection of
sets $S_i\subset \R\times [-H/2,H/2],\, i\in\{1,\ldots,n\}$, satisfying
\begin{align}
\label{eq:Si1}
\bigcup_{\mathbf{k}\in K}\bigcup_{i=1}^n (S_i +\mathbf{k}) = \R\times [-H/2,H/2],
\qquad
\left\vert S_i \cap (S_j + \mathbf{k} ) \right\vert{=0} \; \textrm{for all}\, i\neq j, \, \mathbf{k} \in K,
\end{align}
such that
\begin{align*}
\nabla P(\bx) = \sum_{\mathbf{k}\in K}\sum_{i=1}^n (\z_i+\mathbf{k})\mathds{1}_{S_i+\mathbf{k}}(\bx).
\end{align*}
The corresponding total geostrophic energy is
\begin{align}\label{eq:discreteEnergy1}
	\mathcal{E}&=\frac{f^2}{2}\sum_{\mathbf{k}\in K}\sum_{i=1}^n\int_{\left(S_i+\mathbf{k}\right)\cap\Omega}|\bx- \z_i -\mathbf{k}|^2\, \rd \bx\\
	\label{eq:discreteEnergy2}
	 &\quad -\frac{f^2}{2}\sum_{\mathbf{k}\in K}\sum_{i=1}^n\int_{\left(S_i+\mathbf{k}\right)\cap\Omega}\left(\z_i\cdot \mathbf{e}_2\right)^2\, \rd \bx
	 - \frac{f^2}{2}\int_{\Omega}x_2^2\,\rd \bx
	 + \int_\Omega N^2 \left( x_2 + \frac H2 \right)x_2 \, \rd \bx.
\end{align}
The stability principle says that {at each point in time} this energy is minimised over all periodic rearrangements of particles that preserve potential temperature and absolute momentum. Each such rearrangement corresponds to a collection of sets $S_i\subset \R\times [-H/2,H/2],\, i\in\{1,\ldots,n\}$, satisfying \eqref{eq:Si1} with specified masses
\begin{align}\label{eq:massConstS}
|S_i|=m_i>0 \quad \textrm{for all} \,i\in\{1,\ldots,n\}.
\end{align}
Necessarily,
\begin{align*}
\sum_{i=1}^nm_i=|\Omega|.
\end{align*}
The terms \eqref{eq:discreteEnergy2} are constant over all such collections. This leads to the following minimisation problem:
\begin{align}\label{eq:minProb}
{\min_{\{ S_i\}_{i=1}^n}} \left\{  \sum_{\mathbf{k}\in K}\sum_{i=1}^n\int_{\left(S_i+\mathbf{k}\right)\cap\Omega}|\bx- \z_i -\mathbf{k}|^2\, \rd \bx\, :\, \eqref{eq:Si1}
\,
\text{and}\, \eqref{eq:massConstS}\, \text{hold}\right\}.
\end{align}
We show that the semi-discrete optimal transport problem (Definition \ref{def:SDOT}) is equivalent to \eqref{eq:minProb}. 

Consider $\{S_i\}_{i=1}^n$ satisfying \eqref{eq:Si1} and \eqref{eq:massConstS}. For each $i\in\{1,\ldots,n\}$, define
\begin{align}\label{eq:CiFromSi}
C_i:=\bigcup_{k\in K}\left(S_i+k\right)\cap\Omega.
\end{align}
{Note that $|C_i|=m_i$ for all $i\in\{1,\ldots,n\}$, and $\{C_i\}_{i=1}^n$ is a tessellation of $\Omega$.} Then
\begin{align*}
\sum_{\mathbf{k}\in K}\sum_{i=1}^n\int_{\left(S_i+\mathbf{k}\right)\cap\Omega}|\bx- \z_i -\mathbf{k}|^2\, \rd \bx
&\geq \sum_{\mathbf{k}\in K}\sum_{i=1}^n\int_{\left(S_i+\mathbf{k}\right)\cap\Omega}\min_{\bl\in K}|\bx- \z_i -\bl|^2\, \rd \bx\\
&= \sum_{i=1}^n\int_{C_i}|\bx- \z_i|_{\per}^2\, \rd \bx.
\end{align*} 
By taking the minimum over $\{S_i\}_{i=1}^n$, we see that the minimum \eqref{eq:minProb} is greater than or equal to the minimum attained in the semi-discrete optimal transport problem (Definition \ref{def:SDOT}). 

On the other hand, consider a tessellation $\{C_i\}_{i=1}^n$ of $\Omega$ satisfying the mass constraint $|C_i|=m_i$ for all $i\in\{1,\ldots,n\}$. For each $i\in\{i,\ldots,n\}$ define
\begin{align}\label{eq:SiFromCi}
S_i:=\left\{x\in\bigcup_{\mathbf{k}\in K}\left(C_i+\mathbf{k}\right)\, : \, \underset{\bl\in K}{\mathrm{argmin}} |\bx - {\z_i} - \bl|=\mathbf{0}\right\}.
\end{align}
Elementary calculations show that $\{S_i\}_{i=1}^n$ satisfies \eqref{eq:Si1} and \eqref{eq:massConstS}, and that
\begin{align*}
&C_i = \bigcup_{\bk\in K}(S_i+\bk)\cap \Omega,\\
&\bx\in S_i +\bk \implies |\bx-\z_i|_{\per}=|\bx-\z_i-\bk|.
\end{align*}
Hence
\begin{align*}
\sum_{i=1}^n\int_{C_i}|\bx- \z_i|_{\per}^2\, \rd \bx 
= \sum_{\mathbf{k}\in K}\sum_{i=1}^n\int_{\left(S_i+\mathbf{k}\right)\cap\Omega}|\bx- \z_i|_{\per}^2\, \rd \bx
= \sum_{\mathbf{k}\in K}\sum_{i=1}^n\int_{\left(S_i+\mathbf{k}\right)\cap\Omega}|\bx- \z_i -\bk|^2\, \rd \bx.
\end{align*}
By taking the minimum over $\{C_i\}_{i=1}^n$, we see that the minimum \eqref{eq:minProb} is less than or equal to the minimum attained in the semi-discrete transport problem (Definition \ref{def:SDOT}). 
By combining this with the opposite inequality above, we see that the minimum values are equal, and that the stability principle is equivalent to the semi-discrete optimal transport problem, as claimed. Minimisers are related via \eqref{eq:CiFromSi} and \eqref{eq:SiFromCi}.

\section{Derivation of the ODE}
\label{app:ODEderivation}
{In this appendix we give a formal derivation of the ODE \eqref{eq:ODE1}
from the Lagrangian equation \eqref{eqn:SGESLagCts} by assuming that $P$ is piecewise affine and satisfies the stability principle. 
For a rigorous derivation, see \cite{bourne2022semi}.
By equations \eqref{eq:Z}, \eqref{eqn:pwConstNablaP},
\eqref{eqn:ansatz}, \eqref{eq:brevity}, we have
\begin{align*}
\mathbf{Z}(\mathbf{x}, t) =
\sum_{i=1}^n \left(\z_i(t)	+\mathbf{k}_*(\mathbf{F}(\mathbf{x}, t),\z_i(t))\right)\mathds{1}_{C_{i,\per}(t)}(\mathbf{F}(\mathbf{x}, t)).
\end{align*}
Note that  $\mathbf{k}_*(\mathbf{F}(\mathbf{x}, t),\z_i(t))$ and $\mathds{1}_{C_{i,\per}(t)}(\mathbf{F}(\mathbf{x}, t))$ are piecewise constant in time. Therefore
\[
\partial_t \mathbf{Z}(\mathbf{x}, t)
= \sum_{i=1}^n \dot{\z}_i(t) \mathds{1}_{C_{i,\per}(t)}(\mathbf{F}(\mathbf{x}, t)),
\]
provided that the derivative exists. By \eqref{eqn:cellFlow} and the mass-preserving property of the flow $\mathbf{F}$,
\begin{align}
\label{eq:ODElhs}
\int_{C_{i,\per}(0)}    
\partial_t \mathbf{Z}(\mathbf{x}, t) \, \rd \mathbf{x}
= \int_{C_{i,\per}(t)}    
\sum_{i=1}^n \dot{\z}_i(t) \mathds{1}_{C_{i,\per}(t)}(\mathbf{y})
\, \rd \mathbf{y}
 = \int_{C_{i,\per}(t)}  \dot{\z}_i(t) \, \rd \mathbf{y}
 = \left\vert C_{i,\per}(t) \right\vert \, \dot{\z}_i(t).
\end{align}
On the other hand, integrating the right-hand side of \eqref{eqn:SGESLagCts} over $C_{i,\per}(0)$ gives
\begin{align}
\nonumber
& \int_{C_{i,\per}(0)}    
    J\big(\mathbf{F}(\mathbf{x}, t) - \left(\mathbf{Z}(\mathbf{x}, t) \cdot \mathbf{e}_1\right)\mathbf{e}_1\big) \, \rd \mathbf{x}
\\
    \nonumber
& = \int_{C_{i,\per}(t)}    
    J \mathbf{y} \, \rd \mathbf{y} 
 - \int_{C_{i,\per}(t)}    
    \left[\big( \z_i(t)+\mathbf{k}_*(\mathbf{y},\z_i(t)) \big) \cdot \mathbf{e}_1 \right] J \mathbf{e}_1  \, \rd \mathbf{y}
    \\
    \label{eq:ODErhs}
& = \left\vert C_{i,\per}(t) \right\vert  J \bc_i - 
\left\vert C_{i,\per}(t) \right\vert (\z_i(t) \cdot \mathbf{e}_1) J \mathbf{e}_1,
\end{align}
where $\bc_i$ was defined in \eqref{eq:c}, and 
where we used the fact that
\begin{align*}
  \int_{C_{i,\per}(t)} ( \mathbf{y} - (\mathbf{k}_*(\mathbf{y},\z_i(t)) \cdot \mathbf{e}_1) \mathbf{e}_1 ) \, \rd \mathbf{y} 
 = \int_{C_{i}(t)} \mathbf{y} \, \rd \mathbf{y}.
\end{align*}
Here $C_{i}(t) := C_{i}(\z(t),\w_*(\z(t)))$ is a non-periodic Laguerre cell (see Definition \ref{def:Laguerre}).
By combining \eqref{eqn:SGESLagCts}, \eqref{eq:ODElhs} and \eqref{eq:ODErhs}, we derive \eqref{eq:ODE1} as desired.
}

\section{Linear instability analysis}
\label{Sec:LinearInstabilityAnalysis}
In this section we study the linear instability of the following steady solution (planar Couette flow) of 
equations \eqref{eqn:SGES1} and \eqref{eqn:vThRel0}:
\begin{align}
\label{eq:AA1}
 (u,v,w,\theta,\phi)=
 \left( \overline{u}(z),0,0, \tfrac{N^2 \theta_0}{g}\left(z+\tfrac{H}{2}\right),
\tfrac{N^2}{2}\left(z+\tfrac{H}{2}\right)^2 \right).
\end{align}
This was first done by Eady \cite{eady1949long}, and the unstable perturbations found by Eady were used as initial conditions for the simulations of Williams \cite{williams1967}. We reproduce Eady's results here to make the paper self-contained, to clarify the connection between the 
initial conditions used here and by \cite{williams1967,visram2014framework,yamazaki2017vertical}, and to derive the stable modes, which are not given in \cite{eady1949long}.  

Note that \eqref{eq:AA1} is not the only steady solution of \eqref{eqn:SGES1}, \eqref{eqn:vThRel0}. In fact so is 
\begin{align}
(u,v,w,\theta,\phi)=
\label{eq:SS2}
\left( \overline{u}(z),0,0, \tfrac{\theta_0}{g} \phi'(z),\phi(z) \right)
\end{align}
for any $\phi(z)$.

\subsection{Linear perturbation equations}
Consider the following perturbations of the steady shear flow:
\begin{align}
\label{eq:AA2}
u(x,z,t)& =\overline{u}(z) + \varepsilon u^1(x,z,t), \\
v(x,z,t) & = \varepsilon v^1(x,z,t), \\
w(x,z,t) & = \varepsilon w^1(x,z,t),
\\
\theta(x,z,t) & = \frac{N^2 \theta_0}{g} \left(z + \frac H2 \right) + \varepsilon \theta^1(x,z,t),
\\
\label{eq:AA6}
\phi(x,z,t) & = \frac{N^2}{2} \left( z + \frac H2 \right)^2 + \varepsilon \phi^1(x,z,t).
\end{align}
To linearise equations \eqref{eqn:SGES1} and \eqref{eqn:vThRel0} about the steady solution \eqref{eq:AA1}, we substitute \eqref{eq:AA2}-\eqref{eq:AA6} into \eqref{eqn:SGES1} and \eqref{eqn:vThRel0}, differentiate the resulting equations with respect to $\varepsilon$, and then set $\varepsilon=0$. This results in the following linear PDE for the perturbations 
$(u^1,v^1,w^1,\theta^1,\phi^1)$:
\begin{align}
\label{eq:AA7}
    & \partial_t v^1 + \overline{u} \, \partial_x v^1 + f u^1  = 0,
    \\
    & \partial_t \theta^1 + \overline{u} \, \partial_x \theta^1 + \frac{N^2 \theta_0}{g} w^1 + s v^1 = 0, 
    \\
    & \partial_x u^1 + \partial_z w^1 = 0, 
    \\
    & v^1 = \frac{1}{f} \partial_x \phi^1,
    \\
    \label{eq:AA11}
    & \theta^1  = \frac{\theta_0}{g} \partial_z \phi^1.
\end{align}
This PDE is defined for $(x,z) \in (-\infty,\infty)\times[-H/2,H/2]$. 
We look for solutions that are $2L$-periodic in the $x$-direction and satisfy the rigid-lid boundary condition
\begin{equation}
\label{eq:AA12.5}
    w^1(x,-H/2,t)=w^1(x,H/2,t)=0.
\end{equation}

\subsection{Eigenvalue problem}
By seeking solutions of
\eqref{eq:AA7}-\eqref{eq:AA12.5}
with exponential time dependence $\exp(\omega t)$, $\omega \in \mathbb{C}$, we obtain a generalised eigenvalue problem for the perturbation growth rate $\omega$. Since the eigenvalue problem has $x$-periodic boundary conditions, we can solve it using Fourier series. Consequently, we seek solutions of \eqref{eq:AA7}-\eqref{eq:AA11} of the form
\begin{align}
\label{eq:AA12}
(u^1,v^1,w^1,\theta^1,\phi^1)(x,z,t)
=
(\hat{u}(z),\hat{v}(z),\hat{w}(z),\hat{\theta}(z),\hat{\phi}(z))
E(x,t),
\end{align}
with
\begin{equation*}
E(x,t) = \exp \left( \omega t - \frac{i k \pi x}{L} \right),
\end{equation*}
where $k \in \mathbb{Z}$ is the mode number (and $k \pi/L$ is the wave number). 
Substituting 
\eqref{eq:AA12} into \eqref{eq:AA7}-\eqref{eq:AA12.5} gives the ODE
\begin{align}
\label{eq:AA14}
&  \left( \omega - \frac{i k \pi}{L} \overline{u} \right) \hat{v} + f \hat{u} = 0,
\\
& \left( \omega - \frac{i k \pi}{L} \overline{u} \right)
\hat{\theta} + \frac{N^2 \theta_0}{g} \hat{w} + s \hat{v} = 0,
\\
& - \frac{i k \pi}{L} \hat{u} + \hat{w}' = 0,
\\
\label{eq:AA17}
& \hat{v} = - \frac{1}{f} \frac{i k \pi}{L} \hat{\phi},
\\
\label{eq:AA18}
& \hat{\theta} = \frac{\theta_0}{g} \hat{\phi}',
\end{align}
with boundary conditions
\begin{equation}
\label{eq:AA19}
    \hat{w}(-H/2)=\hat{w}(H/2)=0.
\end{equation}
For $k=0$, $\omega=0$ is an eigenvalue of  \eqref{eq:AA14}-\eqref{eq:AA19} with 
eigenfunctions
$\hat{u}=\hat{v}=\hat{w}=0$, 
$\hat{\phi}$ arbitrary, and
$\hat{\theta}$ determined by \eqref{eq:AA18}.
This eigenvalue corresponds to the family of steady planar Couette flows \eqref{eq:SS2}.

From now on we assume that $k \ne 0$.
We now eliminate $\hat{u}$, $\hat{v}$, $\hat{\phi}$ and $\hat{\theta}$. This can be done by noting that
\begin{align}
\label{eq:elim1}
&\hat{u}=\frac{L}{ik \pi}\hat{w}',
\qquad
\hat{v} = \frac{fL}{ik\pi\left(\frac{ik\pi}{L}\overline{u}-\omega\right)} \hat{w}',
\\ 
\label{eq:elim2}
&\hat{\theta}=\frac{1}{\left(\frac{ik\pi}{L}\overline{u}-\omega\right)}
\left[ \frac{N^2 \theta_0}{g} \hat{w} + \frac{sfL}{ik\pi\left(\frac{ik\pi}{L}\overline{u}-\omega\right)} \hat{w}' \right].
\end{align}
Taking the $z$-derivative of \eqref{eq:AA17}, 
multiplying \eqref{eq:AA18} by $i k \pi g /(\theta_0 f L)$, and adding the resulting equations gives
\begin{equation}
\label{eigen_eady}
f^2 \hat{w}'' + \frac{2 ik \pi f^2 \overline{u}_z }{L \left(\omega - \frac{ik\pi}{L}\overline{u} \right)} \hat{w}'  - \left( \frac{k \pi N}{L} \right)^2 \hat{w} = 0.
\end{equation}
Define the change of variables
\[
W(Z):=\hat{w} \left( \frac{L}{i k \pi \overline{u}_z} \left( \frac{i f \overline{u}_z}{N} Z + \omega \right) \right).
\]
It follows from 
\eqref{eigen_eady} that $W$ satisfies the ODE
\begin{equation}
W''(Z)-\frac{2}{Z}W'(Z)-W(Z)=0.
\end{equation}
The general solution of this ODE is given by
\begin{align}
\label{W(X)}
W(Z) = & A\left(\sinh Z-Z\cosh Z\right)
 +B\left(\cosh Z-Z\sinh Z\right).
\end{align}
This can be used to write the general solution of \eqref{eigen_eady}. First we use the boundary conditions \eqref{eq:AA19} to find the eigenvalues $\omega$. 

Introduce the dimensionless growth rate 
\[
\omega'= \frac{i N\omega}{\bar{u}_zf}  = -\frac{i N \theta_0 \omega}{gs},
\]
and let 
\[
\kappa= \frac{k \pi \mathrm{Bu}}{2} = \frac{k \pi}{2}\frac{NH}{fL}.
\]
Then $\hat{w}(z)=W(\tilde{Z}(z))$, where 
\begin{equation}
    \label{eq:Z(z)}
\tilde{Z}(z) = 
\frac{N}{i f \overline{u}_z} \left(  \frac{i k \pi \overline{u}_z}{L} z - \omega \right)
= \frac{2 \kappa}{H} z + \omega'.
\end{equation}
The boundary conditions \eqref{eq:AA19} 
mean that 
\[
W(\omega' \pm \kappa)=0.
\]
Therefore, from \eqref{W(X)}, we read off that
\begin{align*}
M
\begin{pmatrix}
A \\ B
\end{pmatrix}
=0,
\end{align*}
where $M$ is the 2-by-2 matrix with components
\begin{align*}
M_{11} & = \sinh(\omega'+\kappa) - (\omega'+\kappa) \cosh (\omega'+\kappa), 
\\
M_{12} & = \cosh (\omega'+\kappa) - (\omega'+\kappa) \sinh(\omega'+\kappa), 
\\
M_{21} & = \sinh (\omega'-\kappa) - (\omega'-\kappa) \cosh (\omega'-\kappa),
\\
M_{22} & = \cosh (\omega'-\kappa) - (\omega'-\kappa) \sinh (\omega'-\kappa).
\end{align*}
Non-trivial solutions are then those for which the vector $(A,B)$ belongs to the kernel of $M$.
Imposing the condition that $\det(M)=0$ yields
\begin{equation}
\label{dispersion}
\left(1+\kappa^2-\omp^2\right)\sinh 2\kappa=2\kappa\cosh 2\kappa.
\end{equation}
Since $\kappa\in\mathbb{R}$, any $\omp$ that satisfies the dispersion relation \eqref{dispersion} is either real or purely imaginary. Correspondingly the growth rate $\omega$ is either purely imaginary or real.

\subsection{Unstable modes}
\label{app:unstable}
\begin{figure}
\begin{center}
\includegraphics[scale=1]{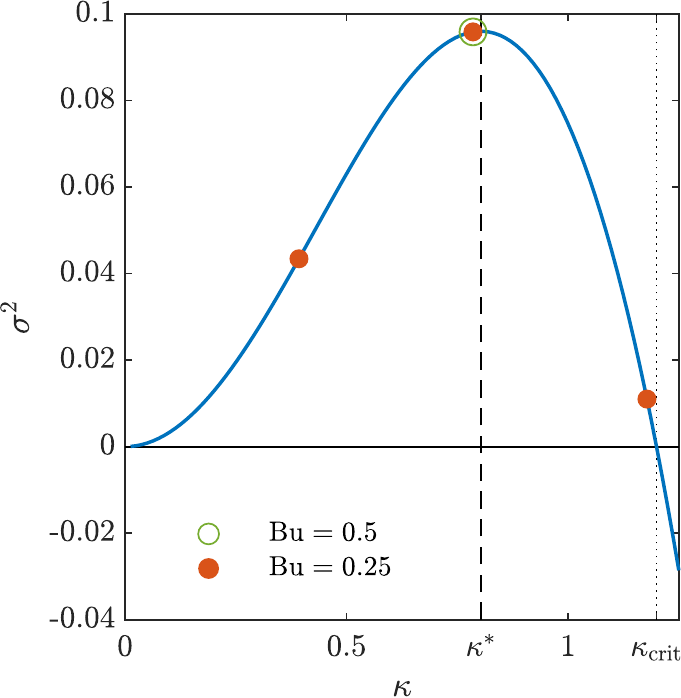}
\end{center}
\caption{
\label{fig:sig} The growth rate $\sigma^2$ from \eqref{eq:sigmaSq}. Only modes $k \in \mathbb{Z}$ with $ |k| \pi \mathrm{Bu}/2<\kappa_{\mathrm{crit}}
=1.19968$
are unstable.
}
\end{figure}

In this section we seek solutions with exponential growth, $\omega > 0$. Set $\omp=i\sigma$ with
\begin{equation}
\label{eq:sigma}
\sigma = -\frac{N \theta_0 \omega}{gs}.
\end{equation}
Note that $\sigma > 0$ since $s<0$. Then  \eqref{dispersion} gives
\begin{align}
\label{eq:sigmaSq}
\sigma^2 & =2\kappa\coth 2\kappa -1-\kappa^2 
= (\kappa - \tanh \kappa)(\coth \kappa - \kappa).
\end{align}
Recall that $\kappa = \kappa_k = k \pi \mathrm{Bu}/2$, where  $k \in \mathbb{Z}$, $k \ne 0$.
Therefore \eqref{eq:sigmaSq} uniquely defines $\sigma>0$ for all $k \in \mathbb{Z}$ such that 
\begin{equation}
\label{eq:admissiblek}    
2\kappa_k\coth2\kappa_k-1-\kappa_k^2 \ge 0.
\end{equation}
(Taking the negative square root $\sigma<0$ in \eqref{eq:sigmaSq} gives an exponentially decaying mode.)
For example, if $\mathrm{Bu}=0.5$, then the only solutions of \eqref{eq:admissiblek} are $k = \pm 1$. If $\mathrm{Bu}=0.25$, then \eqref{eq:admissiblek}  has solutions $k=\pm 1, \pm 2, \pm 3$: see Figure \ref{fig:sig}. 

It is worth observing that there are no exponentially growing modes when 
\[
\mathrm{Bu}>\mathrm{Bu}_{\mathrm{crit}}=\frac{2\kappa_{\mathrm{crit}}}{\pi}
 =0.763739,
\]
where 
$\kappa_{\mathrm{crit}}= 1.19968$ (to 6 s.f.)
is the smallest positive root of \eqref{eq:sigmaSq}. We study this parameter regime in the following section.

It can be shown that  $2\kappa\coth 2\kappa -1-\kappa^2$
is maximised when $\kappa = \pm \kappa^*$ with 
$\kappa^*=0.803058$ (to 6 s.f.).
This is achievable by an integer mode number $k$ if the Burger number satisfies
\begin{align}
\label{eq:BuMax}
    \mathrm{Bu} = \frac{2\kappa^*}{k\pi}
\end{align}
for some $k \in \mathbb{Z}$.
Take $k=1$ and the values of $N$, $L$, $f$, $s$, $\theta_0$, $g$ given in Section \ref{subsec:ParamValues}. We choose $H$ so that \eqref{eq:BuMax} is satisfied. This gives
\[
H = 10224.85 \,\text{m}.
\]
Then the fastest growing mode has growth rate 
\begin{align}
\label{eq:FastestGrowthRate}
    \omega = -\frac{g s}{N \theta_0} \sigma(\kappa^*) = 6.1963 \times 10^{-6} \, \text{s}^{-1}
    = 0.53536 \text{ days}^{-1}.
\end{align}
The fastest growing mode can be rewritten as 
\[
\omega = \frac{\overline{u}_z}{\mathrm{Bu}} \frac{H}{L}  \sigma(\kappa^*).
\]
In particular, we can read off that it is proportional to $\overline{u}_z$ and the aspect ratio $H/L$ of the box, and inversely proportional to $\mathrm{Bu}$.

Next we compute the eigenfunctions for the unstable modes.
By \eqref{eq:Z(z)} we can write
\[
\hat{w}(z)=W(\tilde{z}(z)+i \sigma), \quad \tilde{z}=\frac{2 \kappa}{H}z.
\]
Then
using \eqref{W(X)} we obtain
\begin{align*}
\hat{w}(z) 
=& \tilde{A}\left(\sinh \tilde{z}(z)-\left(\tilde{z}(z)+i\sigma\right)\cosh \tilde{z}(z)\right)
+\tilde{B}\left(\cosh \tilde{z}(z)-\left(\tilde{z}(z)+i\sigma\right)\sinh \tilde{z}(z)\right),
\end{align*}
where $\tilde{A}=A\cos\sigma+iB\sin\sigma$ and $\tilde B=B\cos\sigma+iA\sin\sigma$. 
Applying the boundary conditions \eqref{eq:AA19} gives
\begin{equation*}
\begin{pmatrix}
        S-\left(\kappa+i\sigma\right)C & C-\left(\kappa+i\sigma\right)S \\
        -S+\left(\kappa-i\sigma\right)C & C-\left(\kappa-i\sigma\right)S
    \end{pmatrix}
    \begin{pmatrix}\tilde{A}\\\tilde{B}\end{pmatrix}=
    \begin{pmatrix}0\\0\end{pmatrix},
\end{equation*}
where $S=\sinh\kappa$ and $C=\cosh\kappa$. 
A solution is
\begin{equation*}
\tilde{A}=\frac{g}{N} \sigma,\quad\tilde{B}= \frac{g}{N} i\left( \kappa \coth \kappa - 1\right).
\end{equation*}
(The factor $g/N$ ensures that $\hat{w}$ has the correct units.)
Substituting the expression for $\hat{w}$ into \eqref{eq:elim1}, \eqref{eq:AA17}, \eqref{eq:elim2} gives
\begin{align*}
\hat{u}(z) & = \frac{iN}{f} (\tilde{z}+i\sigma) \left(\tilde{A}\sinh \tilde{z}+\tilde{B}\cosh \tilde{z}\right),
\\
\hat{v}(z) & = -\frac{N^2\theta_0}{gs}\left(\tilde{A}\sinh \tilde{z} +\tilde{B} \cosh \tilde{z} \right),
\\
\hat{\phi}(z) & = 
-\frac{i f L N^2\theta_0}{k \pi g s}\left(\tilde{A}\sinh \tilde{z} +\tilde{B} \cosh \tilde{z} \right),
\\
\hat{\theta}(z) & = 
-\frac{iN^3\theta_0^2}{g^2 s}\left(\tilde{A}\cosh \tilde{z}+\tilde{B}\sinh \tilde{z} \right).
\end{align*}
Therefore the unstable perturbations are 
\begin{gather*}
u^1(x,z,t)=\Re\left[\hat{u}(z)\exp \left(\omega t - \frac{i k \pi x}{L}\right)\right], \quad \textrm{etc}
\end{gather*}
which yields
\begin{align*}
& u^1 = 
e^{\omega t} \frac{g}{f} \left[   
\left(  - A_2^2 \sinh \tilde z -A_1 \tilde z \cosh \tilde z \right) \cos
\tfrac{k \pi x}{L} +   
\left( A_2 \tilde z \sinh \tilde z - A_1 A_2 \cosh \tilde z \right) \sin \tfrac{k \pi x}{L}
\right],
\\
& v^1  = e^{\omega t} \frac{N \theta_0}{s} \left[-A_2 \sinh \tilde z \cos \tfrac{k \pi x}{L}  - A_1 \cosh \tilde z \sin \tfrac{k \pi x}{L} \right],
\\
& w^1  =
e^{\omega t} \frac{g}{N} \left[ 
-A_2 \cosh \tilde z \left( \tilde z \cos \tfrac{k \pi x}{L}
+ A_2 \sin \tfrac{k \pi x}{L}
\right) \right.
+ \kappa \coth \kappa \sinh \tilde z \cos \tfrac{k \pi x}{L}
\\
& \qquad \qquad \qquad \qquad \qquad \qquad \qquad \qquad \qquad \qquad \qquad
+ \left. A_1 \left( \cosh \tilde z - \tilde z \sinh \tilde z \right) \sin 
\tfrac{k \pi x}{L} \right],
\\
& \theta^1  =
e^{\omega t} \frac{N^2 \theta_0^2}{g s} 
\left[ -A_2 \cosh \tilde z \sin \tfrac{k \pi x}{L} + A_1 \sinh \tilde z \cos \tfrac{k \pi x}{L} \right],
\\
& \phi^1  =
e^{\omega t} \frac{fLN \theta_0}{k \pi s} 
\left[ -A_2 \sinh \tilde z \sin \tfrac{k \pi x}{L} + A_1 \cosh \tilde z \cos \tfrac{k \pi x}{L} \right],
\end{align*}
where
\[
A_1 = \kappa \coth \kappa - 1, \qquad A_2 = \sigma.
\]
These expressions agree with those in \cite{williams1967}.

In Section \ref{Subsec:UNM} we use the following initial condition for the gradient of the modified geopotential:
\begin{align}
\label{eq:gradPinitial}
\nabla P((x,z),0) = 
\begin{pmatrix}
\frac{1}{f} v^0 +x
\\
\frac{g}{f^2 \theta_0} \theta^0(z) 
\end{pmatrix}
+ \lambda\begin{pmatrix}
\frac{1}{f} v^1(x,z,0)
\\
\frac{g}{f^2 \theta_0} \theta^1(x,z,0)
\end{pmatrix},
\end{align}
where $v^0=0$ and $\theta^0(z)=\frac{N^2\theta_0}{g} \left( z + \frac H2 \right)$ are the steady meridional velocity and potential temperature, $k=1$, and 
\begin{equation}
\label{eq:lambda}
\lambda =
\frac{sa}{N \theta_0}.
\end{equation}
In particular, $v_{\mathrm{u}}(x,z)=\lambda v^1(x,z,0)$ and $\theta_{\mathrm{u}}(x,z)=\lambda \theta^1(x,z,0)$.

\subsection{Stable modes}
\label{app:stable}
In this section we study the parameter regime $\mathrm{Bu}>\mathrm{Bu}_{\mathrm{crit}}$, where the eigenvalues $\omega$ are purely imaginary for all $k \in \mathbb{Z}$. This means that the steady Couette solution \eqref{eq:AA1} is 
linearly stable with respect to normal mode perturbations.

It turns out, however, that the normal modes do not form a complete basis, in the following sense. For each mode number $k$, there are only two solutions $\pm \sigma(\kappa_k)$ of \eqref{eq:sigmaSq}, and therefore only two eigenvalues $\pm \omega_k$. Consequently the boundary value problem \eqref{eigen_eady}, \eqref{eq:AA19} only has two eigenvalues, and not a countable set of eigenvalues as one might expect. This means that it is not possible to represent every initial condition of the linear perturbation equations \eqref{eq:AA7}-\eqref{eq:AA11} as a sum of normal modes (eigenfunctions). Therefore stability with respect to normal mode perturbations does not guarantee stability with respect to all perturbations. The origin of this problem is the denominator $\omega - \tfrac{i k \pi}{L} \overline{u}$ in equation \eqref{eigen_eady}, which can vanish.

For a closely related but simpler baroclinic instability problem, Pedlosky \cite{pedlosky1964initial} showed that any initial perturbation can be represented by supplementing the discrete spectrum of normal modes with a continuous spectrum. By doing this he showed that the linear stability of the steady solution is in fact determined by the incomplete set of normal modes. Our simulations suggest that the same is true for our problem. We plan to explore this in a future paper.

To compute the normal modes in the parameter regime $\mathrm{Bu}>\mathrm{Bu}_{\mathrm{crit}}$, we write $\sigma = i \gamma$, $\gamma \in \mathbb{R}$.
Since $\mathrm{Bu}>\mathrm{Bu}_{\mathrm{crit}}$, then 
$|\mathrm{\kappa}|>\mathrm{\kappa}_{\mathrm{crit}}$ for all $k \in \mathbb{Z}$, $k \ne 0$, and hence the right-hand side of  \eqref{eq:sigmaSq} is negative. Therefore we read off from equation \eqref{eq:sigmaSq} that
\[
\gamma^2 = |(\kappa - \tanh \kappa)(\coth \kappa - \kappa)|.
\]
Then, by \eqref{eq:sigma},
\[
\omega = \tilde{\omega} i, \qquad \tilde{\omega} = - \frac{gs \gamma}{N \theta_0} \in \mathbb{R}.
\]
It follows that 
\begin{align*}
v^1(x,z,t) =\Re\left[\hat{v}(z)e^{i \left(\tilde{\omega} t - \frac{k \pi x}{L}\right)}\right], 
\qquad
\theta^1(x,z,t) =\Re\left[\hat{\theta}(z)e^{i \left(\tilde{\omega} t - \frac{k \pi x}{L}\right)}\right],
\end{align*}
which yields
\begin{align*}
& v^1  = 
\frac{N \theta_0}{s} ( \gamma \sinh \tilde z + A_1 \cosh \tilde z) \sin \left( \tilde \omega t - \tfrac{k \pi x}{L} \right),
\\
& \theta^1  = \frac{N^2 \theta_0^2}{gs} 
( \gamma \cosh \tilde z + A_1 \sinh \tilde z) \cos \left( \tilde \omega t - \tfrac{k \pi x}{L} \right).
\end{align*}
Substituting these expressions into \eqref{eq:gradPinitial} and taking $k=1$ gives the initial condition that we use in Section \ref{Subsec:SNM}. 
{In particular, $v_{\mathrm{s}}(x,z)=\lambda v^1(x,z,0)$ and $\theta_{\mathrm{s}}(x,z)=\lambda \theta^1(x,z,0)$, where $\lambda$ was defined in \eqref{eq:lambda}.}
Observe that this perturbation is a travelling wave with wave speed
\[
c_{k} = \frac{\tilde \omega}{\tfrac{ k \pi }{L}} = - \frac{g s \gamma L}{k \pi N \theta_0}.
\]
The limiting value of the wave speed for small wavelength perturbations is
\begin{align}
\label{eq:stablewavespeed}
c_{\infty} & = \lim_{k \to \infty} c_k
=
- \frac{g s L}{\pi N \theta_0} \lim_{k \to \infty} \left( \left| \frac{\kappa - \tanh \kappa}{k} \right|  \left| \frac{\coth \kappa - \kappa}{k} \right| \right)^{\frac 12} 
= 
- \frac{g s L}{\pi N \theta_0} \frac{\pi \mathrm{Bu}}{2}
= - \frac{gsH}{2 f \theta_0}.
\end{align}

\section{Derivatives of $\mathcal{K}$}\label{app}
{Recall that $\mathcal{K}$ was defined in equation \eqref{eqn:KantFunct}.}
The gradient and Hessian of $\mathcal{K}$ are given by
\begin{align*}
	\nabla\mathcal{K}&=\obm-\m ,\\ D^2\mathcal{K}&=-D_{\w}\m=-\left(\pdonet{m_i}{w_j}\right)_{i,j=1}^n,
\end{align*}
where $\m$ is the cell-area map from Definition \ref{defn:cellAreas}. 
{For a proof, see for example \cite{kitagawa2016convergence,merigot2021optimal}.}
The system \eqref{eqn:NewtDirAlg} defining the $k$\textsuperscript{th} Newton direction $\mathbf{d}^{(k)}$ can therefore be written as
\begin{align*}
	\left\{\begin{array}{l}
		- \mathrm{D}_{\w} \m\left(\z,\w^{(k)}\right) \mathbf{d}^{(k)}=\m\left(\z,\w^{(k)}\right)-\obm,\vspace{3pt}\\
		\mathbf{d}^{(k)}\cdot \mathbf{e}_n=0.
	\end{array}\right.
\end{align*}
{For $i \ne j$, 
\begin{align*}
\frac{\partial m_i}{\partial w_j}
= -\frac 12 \sum_{\mathbf{k}\in K} 
\frac{\mathrm{length}(e(i,j,\mathbf{k}))}{|\z_i + \mathbf{k} - \z_j|},
\end{align*}
where $e(i,j,\mathbf{k})$
is the edge between the non-periodic Laguerre cells with seeds $\z_i+\mathbf{k}$ and $\z_j$, which is defined by
\begin{multline*}
e(i,j,\mathbf{k})
=
\Big\{ \mathbf{x} \in \mathbb{R} \times \left[-\tfrac H2, \tfrac H2\right] : 
\forall \; m \in \{1,\ldots,n\}, \; \mathbf{l} \in K,
\\
| \mathbf{x}-(\z_i + \mathbf{k})|^2 - w_i \le 
| \mathbf{x}-(\z_m+\mathbf{l})|^2 - w_m, 
\; | \mathbf{x}-\z_j |^2 - w_j \le 
| \mathbf{x}-(\z_m+\mathbf{l})|^2 - w_m \Big\}.
\end{multline*}
Note that this set may be empty, in which case $\partial m_i/\partial w_j = 0$.
The diagonal entries of $\mathrm{D}_{\w} \m$ are
\begin{align*}
\frac{\partial m_i}{\partial w_i}
& = -\sum_{\substack{j=1 \\ j \ne i}}^n \frac{\partial m_i}{\partial w_j}.
\end{align*}
The expressions for the derivatives of $\m$ are proved for example in \cite{kitagawa2016convergence,merigot2021optimal} 
for non-periodic Laguerre tessellations and in \cite{BournePearceRoper} for periodic Laguerre tessellations.
}

\section*{Acknowledgements}
C. P. Egan is supported by The Maxwell Institute Graduate School in Analysis and its Applications, a Centre for Doctoral Training funded by the UK Engineering and Physical Sciences Research Council (EPSRC) via the grant EP/L016508/01, the Scottish Funding Council, Heriot-Watt University and the University of Edinburgh. 
D. P. Bourne is supported by the EPSRC via the grant EP/V00204X/1 Mathematical Theory of Polycrystalline Materials.
C. J. Cotter would like to acknowledge the NERC grants NE/K012533/1 and NE/M013634/1.
B. Pelloni and M. Wilkinson gratefully acknowledge the support of the EPSRC via the grant EP/P011543/1 Analysis of models for large-scale geophysical flows.

\bibliographystyle{siam}
\bibliography{SGEadySliceBib}

\end{document}